\documentclass{amsart}
\usepackage{amsmath}
\usepackage{amssymb}
\usepackage{amsthm}
\usepackage{graphicx}
\usepackage{mathtools}
\usepackage{multirow}

 \usepackage[dvipsnames]{xcolor}

 \usepackage{color}

\newcommand{\margincolor}{red}      
\definecolor{darkgreen}{rgb}{0,0.7,0}

\newcounter{margincounter}
\newcommand{\marginnum}{
\ifnum\value{margincounter}<10
\textcolor{\margincolor}{\begin{picture}(0,0)\put(2.2,2.4){\circle{9}}\end{picture}\footnotesize\arabic{margincounter}}
\else\ifnum\value{margincounter}<100
\textcolor{\margincolor}{\begin{picture}(0,0)\put(4.256,2.5){\circle{11}}\end{picture}\footnotesize\arabic{margincounter}}
\else
\textcolor{\margincolor}{\begin{picture}(0,0)\put(6.8,2.5){\circle{14}}\end{picture}\footnotesize\arabic{margincounter}}
\fi\fi
}

\usepackage{hyperref}

 \newtheorem{thm}{Theorem}[section]
 \newtheorem{cor}[thm]{Corollary}
 \newtheorem{lem}[thm]{Lemma}
 
 \newtheorem{conj}[thm]{Conjecture}

 \theoremstyle{definition}
 \newtheorem{defn}[thm]{Definition}
 
 \theoremstyle{remark}
 \newtheorem{rem}[thm]{Remark}
 
 \numberwithin{equation}{section}

 \newcommand{\R}{\mathtt{r}}
 \newcommand{\T}{\mathtt{t}}
 
 \newcommand{\B}{\mathtt{b}}
 \newcommand{\C}{\mathtt{c}}
 \newcommand{\s}{\mathtt{s}}

\title{Universal Geometric Coefficients \break for the Dread Torus}
\author{Jessica N. Copher}

\begin{document}

\begin{abstract} 
The combinatorial structure of certain cluster algebras of geometric type can be modeled by \emph{marked surfaces}. 
We show that a particular marked surface (the \emph{dread torus}) has the Null Tangle Property. 
This result implies that universal geometric coefficients (over $\mathbb{Z}$ and $\mathbb{Q}$) exist for every initial exchange matrix associated with the dread torus and allows us to explicitly construct these coefficients. 
To put this result in perspective:  Only relatively few marked surfaces are known to have the Null Tangle Property, and, besides these, only three others are known to have universal geometric coefficients. 
We also describe the dread torus' combinatorics by characterizing its \emph{arcs}, \emph{triangulations}, and \emph{allowable curves} (in preparation for the proof that the dread torus has the Null Tangle Property) and describing its \emph{mutation fan}. 
\end{abstract}

\thanks{This paper is a shortened version of a chapter in the author's Ph.D. dissertation \cite{mydissertation}. The author wishes to thank her Ph.D. advisor, Nathan Reading, for his numerous editorial suggestions and his sound mathematical guidance regarding this paper. She also wishes to thank him for helping generate the mutation fan figures and helping create final, camera-ready versions of a few other figures.}

\maketitle

\tableofcontents

\section{Introduction} \label{Sect-Intro}

A \emph{cluster algebra} is a type of ring that has certain generators called \emph{cluster variables} grouped into same-sized sets called \emph{clusters} \cite{CA1}. 
The clusters are determined by a recursive process \cite{CA1, CA4}. 
First, choose an initial cluster $x_1, x_2, ..., x_n$ and some initial ``combinatorial data'' consisting of an $n \times n$ skew-symmetrizable matrix (called the \emph{initial exchange matrix}) and $n$ \emph{initial coefficients} and then use that data to perform \emph{mutations} to that cluster-data pair. 
(Each mutation 
$\mu_k$ 
is performed in three steps. 
The first step, called \emph{matrix mutation}, uses the old exchange matrix to construct a new one. 
The other two steps yield $n$ new coefficients and a new cluster.)  
Continue performing mutations to each new cluster-data pair. 

A cluster algebra is said to be \emph{geometric} or \emph{of geometric type} (in the broad sense introduced in \cite{26}) if its coefficients are \emph{geometric}, i.e., if they come from a \emph{tropical semifield over $R$} (where $R$ is usually $\mathbb{Z}$, $\mathbb{Q}$, or $\mathbb{R}$). 
The initial coefficients $y_1, y_2, \dots, y_n$ of a geometric cluster algebra (with initial exchange matrix $B$) are called \emph{universal (over $R$)} if the coefficients of any geometric cluster algebra that has $B$ as its initial exchange matrix (and whose coefficients come from a tropical semifield over $R$) can be obtained from $y_1, y_2, \dots, y_n$ via a unique \emph{coefficient specialization}. 
(See \cite{CA4, 26} for more details.) 

Certain geometric cluster algebras can be modeled by \emph{marked surfaces} \cite{CATS1, CATS2}, each of which consists of a compact, oriented surface (with or without boundary) and a nonzero finite number of points on the surface that are designated \emph{marked points}. 
If a marked point occurs in the interior of a surface, it is called a \emph{puncture}.
In this model, an initial exchange matrix is specified by fixing a marked surface and a \emph{(tagged) triangulation} of that surface 
and initial geometric coefficients are specified by fixing a collection of \emph{laminations} of that surface \cite{CATS1, CATS2}.

The focus of this paper is a particular surface $S$, namely, a torus with a removed open disk $D$, which has one marked point $p$ on its boundary.  We denote this marked surface by $(S,p)$ and refer to it as the \emph{dread torus}.\footnote{``Dread torus'' is a shortening of the name ``dreaded torus,'' which was coined extemporaneously in a talk by Gregg Musiker (in reference to a difficult case of a proof) and, since then, has been used occasionally in the literature (see \cite{Bucher-etal}, \cite{Canakci-Lee-Schiffler}, \cite{Lee-Lee-Mills}, \cite{Matherne-Muller}). We prefer this shortened version since, according to \cite{Merriam-Webster}, ``dread'' has a positive interpretation (i.e., ``inspiring awe'') whereas ``dreaded'' has only negative interpretations
(e.g., inspiring ``extreme reluctance to meet or face'').
}
The main result of this paper is that $(S,p)$ has the Null Tangle Property (Theorem~\ref{Thm-NTP}), which implies that universal geometric coefficients exist for any initial exchange matrix associated with a triangulation of $(S,p)$ and can be explicitly identified (as we do in Corollary~\ref{Cor-universal-geom-coeffs-wrt-Delta0} for one triangulation $\Delta_0$ of $(S,p)$). 
(The universal geometric coefficients associated with the other triangulations of $(S,p)$ can be easily derived from Corollary~\ref{Cor-universal-geom-coeffs-wrt-Delta0} by matrix mutation; see Remark~\ref{Rem-universal-geom-coeffs-wrt-notDelta0}.) 

To put this result in perspective, here is a list\footnote{This list corrects an error in \cite{27}, namely, mistakenly classifying the unpunctured quadrilateral and once-punctured digon as having the Null Tangle Property rather than as being null surfaces.} of all the marked surfaces that, prior to this work, were known to have the Null Tangle Property:  
\begin{enumerate}
	\item \label{Thm27-7.4a} an unpunctured disk with at least five marked points on its boundary,
	\item \label{Thm27-7.4b} a once-punctured disk with at least three marked points on its boundary,
	\item \label{Thm27-7.4c} an unpunctured annulus,
	\item \label{Thm27-7.4d} a twice-punctured disk,
	\item \label{Thm28-3.2} a once-punctured torus,
	\item \label{Thm34-1.1} a four-punctured sphere,
        \item \label{Thm27-7.4e} a once-punctured annulus, and
        \item \label{Thm27-7.4f} a sphere with three removed open disks
\end{enumerate}
(see \cite[Theorem 1.1]{34} for (\ref{Thm27-7.4a})-(\ref{Thm27-7.4d}), (\ref{Thm27-7.4e}), and (\ref{Thm27-7.4f}); \cite[Theorem 7.4]{27} for (\ref{Thm28-3.2}); and \cite[Theorem 3.2]{28} for (\ref{Thm34-1.1})). 
Only three marked surfaces (collectively called the \emph{null surfaces}) are known (to date) to \emph{not} have the Null Tangle Property: the once-punctured monogon, the once-punctured digon, and the unpunctured quadrilateral \cite[Subsection 1.1]{45}. 
The initial exchange matrices associated with these three surfaces have universal geometric coefficients for trivial reasons (see \cite[Theorem~2.1]{27} and \cite[Theorem~3.13 and Remark~3.14]{45}). 

To prove that $(S,p)$ has the Null Tangle Property, we first characterize the \emph{arcs} and triangulations of $(S,p)$ (which correspond to the cluster variables and clusters, respectively) and then characterize the \emph{allowable curves} of $(S,p)$ (which are the building blocks for laminations and, hence, for geometric coefficients). 
After proving that $(S,p)$ has the Null Tangle Property, we describe the \emph{mutation fan} for $\Delta_0$ (i.e., a geometric encoding of the mutations of the initial exchange matrix associated with $\Delta_0$). 
These secondary results, along with the universal geometric coefficients given in Theorem~\ref {Thm-shear-coords-of-ACs}, provide a fairly complete combinatorial picture of $(S,p)$. 

To put these secondary results in perspective, thorough combinatorial descriptions exist for the first six of the eight surfaces (listed above) that are known to have the Null Tangle Property.  
(The first four surfaces have the combinatorial structure of generalized associahedra of types $A$, $D$, $\widetilde{A}$, and $\widetilde{D}$, respectively \cite{CATS1, CA2, 40, 47}; the fifth and sixth surfaces have their combinatorial structures described, in \cite{28} and \cite{34}, respectively, in manners similar to the manner of this paper.) 
The combinatorial structures of the seventh and eighth surfaces are less well understood. 

\section{Arcs and Triangulations}

In this section, we characterize the dread torus' arcs and triangulations. 
This is accomplished by utilizing the close relationship between the dread torus and the \emph {once-punctured torus} (i.e., a torus $T$ with a single marked point).  We denote the once-punctured torus by $(T,p)$ (this means that ``$p$'' will sometimes denote the once-punctured torus' marked point and will sometimes denote the dread torus' marked point, but it should always be clear from context which one is meant). 

Some of the definitions given here are specific to $(S,p)$ and $(T,p)$, rather than applying to marked surfaces in general. 
The arcs of $(T,p)$ are closely related to Farey points. 
We follow the treatment from \cite[Section 4]{28}.

\begin{defn} \label{Defn-Farey-stuff}
A \emph{Farey point} is an integer point $(a,b)$ such that $a$ and $b$ are relatively prime (i.e., $a$ and $b$ are not both $0$ and $\gcd(a, b) = 1$).  A \emph{standard Farey point} is a Farey point such that either $a > 0$ or $a=0$ and $b=1$.
The \emph{standard form} of a rational or infinite slope is a fraction $b/a$ such that $(a, b)$ is a standard Farey point. 
In this paper, all fractions are written in standard form unless otherwise specified.
\end{defn}

\begin{defn} [Cf. \protect{\cite[Definition 2.2]{CATS1}}]
An \emph{arc} in either $(S, p)$ or $(T, p)$ is a curve whose endpoints are both at the marked point $p$ but that otherwise lies in the interior of the surface and is non-self-intersecting; however, any curve that forms the boundary of a monogon (i.e., is contractible) or that, together with the boundary of $S$, forms a digon (i.e., is contractible to the boundary of $S$) is not an arc.  Arcs in $(S, p)$ and $(T, p)$ are considered up to isotopy fixing $p$.
\end{defn}

Recall that the torus can be realized as the quotient space $\mathbb{R}^2/\sim$ where $\sim$ is the equivalence relation defined by $(x,y) \sim (x',y')$ if and only if $x-x', y-y' \in \mathbb{Z}$.  
Moreover, $\mathbb{R}^2$ is the universal cover for the torus with the natural map (the map that sends $(x,y)$ to its equivalence class) as its covering map.  The deck transformations of this covering map are the integer translations of $\mathbb{R}^2$.   
We will take the image of the integer points under this covering map to be $p$ in the once-punctured torus.   

The following theorem characterizes the arcs in $(T, p)$.  It is part of \cite[Proposition 4.1]{28}, recast in terms of standard Farey points.

\begin{thm} \label{Prop-4-1-from-28}
The arcs in $(T, p)$ are the images, under the universal covering map from $\mathbb{R}^2$ to $T$, of the line segments connecting the origin to the standard Farey points.  Moreover, the map from line segments to arcs is bijective (up to isotopy).
\end{thm}

For each arc in the once-punctured torus, we refer to its lift to $\mathbb{R}^2$ that has the origin and a standard Farey point as endpoints as its \emph{standard lift}. 
A standard lift has slope $b/a$ if and only if its non-origin endpoint is $(a, b)$. 
For simplicity, in this case, we say that the arc has slope $b/a$ as well. 
We always orient standard lifts as heading from the origin to the Farey point. 

Let $\mathbb{R}^2\setminus D\mathbb{Z}^2$ denote the real plane with identical small open disks removed at each integer point so that the boundary of each removed disk is tangent to the vertical line passing through the integer point and lies on that line's left. 
The \emph{boundary components} of $\mathbb{R}^2 \setminus D\mathbb{Z}^2$ are the boundaries of the removed disks. 
Observe that $\mathbb{R}^2\setminus D\mathbb{Z}^2$ is a (non-universal) cover for $S$ under the map analogous to the universal covering map from $\mathbb{R}^2$ to $T$.  We will take $p$ to be the boundary point of~$S$ that is the image of the integer points under this covering map.

The following theorem characterizes the arcs in $(S, p)$. The essence of its proof was previously given as part of a proof of a different theorem in \cite[Subsection~4.3]{30}, including the key tactic of thinking of the dread torus' boundary as being so small that it can be considered a ``thickened point.''  

\begin{thm} \label{Thm-arcs-as-images}
The arcs in $(S, p)$ are the images, under the covering map, of the curves in $\mathbb{R}^2\setminus D\mathbb{Z}^2$ that start at the origin, end at a standard Farey point $(a, b)$, and have the following three components:

\begin{enumerate}
    \item \label{origin-spirals} $k$ clockwise spirals around the boundary component at the origin, starting at the origin and ending on the line segment between the origin and $(a, b)$ (partial spirals are rounded to the nearest whole number, $k=0$ indicates no spirals, and negative values of $k$ indicate counterclockwise spirals).
    \item \label{non-origin-spirals} a set of spirals around the boundary component at $(a, b)$, ending at $(a,b)$, that does not intersect the translate of the curve formed by shifting the curve $a$ units to the right and $b$ units up.  
    \item \label{segment} the line segment of slope $b/a$ that joins the spirals described in (\ref{origin-spirals}) and (\ref{non-origin-spirals}). 
\end{enumerate}
This map from curves to arcs is bijective (up to isotopy). 
\end{thm}

The line segment in (\ref{segment}) is only defined if the size of the removed disks of $\mathbb{R}^2\setminus D\mathbb{Z}^2$ is sufficiently small. 
Thus, in practice, we either fix a size for the removed disks that is sufficiently small for all arcs under consideration (this only works if the arcs correspond to a finite number of slopes) or else permit the line segments to bend slightly around the components described in (\ref{origin-spirals}) and (\ref{non-origin-spirals}) (this is necessary if the arcs correspond to an infinitude of slopes). 
In the theorem's proof, when we think of the boundary of $S$ as being a ``thickened point,'' we are, in essence, fixing a size for the removed disks that is sufficiently small for \emph{all} the arcs of $(S,p)$. 

Immediately following this theorem's proof, we will show that, for fixed~$k$, there are exactly two possible sets of spirals around the boundary component at $(a,b)$.  

\begin{proof} 
A lift of an arc in $(S,p)$ is precisely a non-contractible curve in $\mathbb{R}^2 \setminus D \mathbb{Z}^2$ that does not, along with a boundary component, form a digon and whose endpoints are both integer points but that otherwise lies in the interior of $\mathbb{R}^2 \setminus D \mathbb{Z}^2$ and does not intersect itself or any other lift of the same arc except possibly at its endpoints.

If we think of the dread torus' boundary as being a ``thickened point,'' the dread torus becomes, essentially, the once-punctured torus. 
Thus, 
the arcs of the dread torus are identical to those of the once-punctured torus\textemdash except, possibly, near the boundary.  
By Theorem \ref{Prop-4-1-from-28}, the arcs of the once-punctured torus are the images, under the universal covering map from $\mathbb{R}^2$ to $T$, of the line segments connecting the origin to the standard Farey points with exactly one such line segment per arc.  
Thus, each arc of the dread torus is the image, under the covering map from $\mathbb{R}^2 \setminus D\mathbb{Z}^2$ to $S$, of a curve that connects the origin to a standard Farey point and coincides with a straight line segment between these points everywhere except, possibly, within small open neighborhoods $U$ and $V$ around the boundary components at the origin and at the standard Farey point, respectively.

Let $(a,b)$ be a standard Farey point, let $P$ be the intersection point of the boundary of $U$ and the line segment joining the origin and $(a,b)$, and let $Q$ be the intersection point of the boundary of $V$ and the line segment joining the origin and $(a,b)$.
Inside $U$, it is evident that there are many non-isotopic ways that a non-self-intersecting curve in the interior of $\mathbb{R}^2 \setminus D \mathbb{Z}^2$ can join the origin to the point $P$, namely, the curve can complete any integer number $k$ of spirals around the boundary component at the origin before leaving $U$.  Likewise, inside $V$, a non-self-intersecting curve in the interior of $\mathbb{R}^2 \setminus D \mathbb{Z}^2$ can join 
the point $(a,b)$ to the point $Q$ via any integer number of spirals around the boundary component at the point $(a,b)$.  
However, since the deck transformations of the covering map of the dread torus are precisely the integer translations, the lifts of any given arc consist of all the integer translations of the lift whose endpoints are the origin and $(a,b)$.  
Thus, the only curves from the origin to $(a,b)$ that are lifts of arcs in the dread torus are the ones that consist of (a) $k$ spirals around the boundary component at the origin, (b) the line segment from $P$ to $Q$, and (c) an integer number of spirals around the boundary component at $(a,b)$ such that the spirals do not intersect the translation of the curve that starts at $(a,b)$.  Moreover, all such curves are lifts of arcs of the dread torus. 
\end{proof}

As with the once-punctured torus, we refer to the lift of each arc in the dread torus that has the origin and a standard Farey point as its endpoints as its \emph{standard lift} and always orient these lifts as heading from the origin to the Farey point.

\begin{figure}
\scalebox{0.3}{\includegraphics{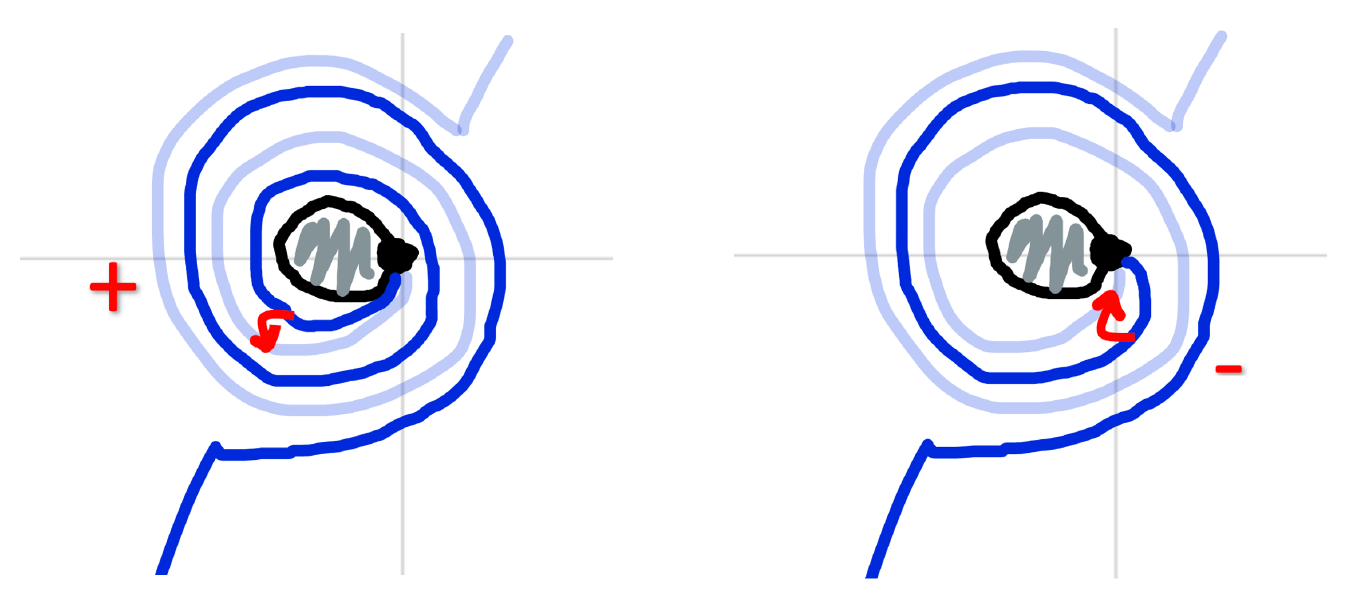}}
\caption{$C$ (light blue) and $C'$ (dark blue) near the origin in $\mathbb{R}^2 \setminus D\mathbb{Z}^2$.}
\label{arc-ends}
\end{figure}

To fully characterize the possible end behavior of arcs in $(S, p)$, let $C$ denote the standard lift of an arc in $(S,p)$ and let $C'$ denote the translate of $C$ that ends at the origin. 
Notice that Condition \ref{non-origin-spirals} of Theorem~\ref{Thm-arcs-as-images} is equivalent to the condition that the ending set of spirals of $C'$ does not intersect the beginning set of spirals of~$C$.  
As shown in Figure~\ref{arc-ends}, there are only two ways this equivalent condition can be met: either the beginning of $C$ lies counterclockwise (CCW) of the end of $C'$ in~$S$ or else the beginning of $C$ lies clockwise (CW) of the end of $C'$ in $S$. 

\begin{defn} \label{defn-sign}
We refer to the end behavior of a lift of an arc in $(S,p)$ as the lift's \emph{sign}.  Lifts with the first kind of end behavior described above have a sign of~$+$ (see the left picture of Figure~\ref{arc-ends}), and lifts with the second kind of end behavior described above have a sign of~$-$ (see the right picture of Figure~\ref{arc-ends}).  We use the order $ - < +$ on $\{-, +\}$.
\end{defn}

As a result of Theorem~\ref{Thm-arcs-as-images} and the above discussion, each arc in the dread torus can be uniquely identified by listing three attributes of its standard lift:
(a) the slope $b/a$ of the line joining the two sets of spirals (which specifies a standard Farey point), 
(b) the number $k$ of spirals around the boundary component at the origin, and
(c) the sign $\sigma$.
We say, for example, ``the standard lift $(b/a, k, \sigma)$ of the arc'' or ``the arc given by $(b/a, k, \sigma)$.''  For simplicity, we also say that the arc has slope $b/a$, spiral number $k$, and sign $\sigma$. 
We use the following total order on arcs of the same slope:  If $\gamma_1$ and $\gamma_2$ are given, respectively, by $(m, k_1, \sigma_1)$ and $(m, k_2, \sigma_2)$, we   write ``$\gamma_1 < \gamma_2$'' (and say, ``$\gamma_1$ is less than $\gamma_2$'') if either (a) $k_1 < k_2$ or (b) $k_1 = k_2$ and $\sigma_1 < \sigma_2$. 

\begin{defn}[Cf. \protect{\cite[Definition 2.4]{CATS1}}]
Two arcs in either $(S,p)$ or $(T,p)$ are \emph{compatible} if they have isotopy representatives that do not intersect except at $p$.
\end{defn}

The next two theorems, together, characterize the pairs of compatible arcs in $(S,p)$. (We will use these theorems in characterizing the triangulations.) 
To prove the first theorem, we make use of the following lemma. 
The proof of this lemma is straightforward and can be found in \cite{mydissertation}. 

\begin{lem} \label{Lem-change-spiral-numbers}
For every integer $n$, there is a homeomorphism $f_n$ of $S$ that fixes the boundary of $S$ pointwise and increases the number of spirals of every arc by $n$ but does not, up to isotopy, otherwise alter the arcs.
\end{lem}

The map $f_n$ described in Lemma~\ref{Lem-change-spiral-numbers} is an instance of a \emph{Dehn twist} (as defined in \cite{Farb-Margalit}).

\begin{thm}  \label{Thm-compatible-arcs-same}
Two arcs in $(S,p)$ given by $C_1=(m, k_1, \sigma_1)$ and $C_2=(m, k_2, \sigma_2)$ with $\sigma_1 \geq \sigma_2$ are compatible if and only if $\sigma_1 > \sigma_2$ and $k_2 \in \{k_1, k_1+1 \}$ (or, in other words, if and only if there is a cover relation between them in the total order on arcs of slope $m$).
\end{thm}

To explain the proof of this theorem and the next, we rely heavily on Figures~\ref{Fig2a-2b-2f} and \ref{Fig1a-1b}. 
Each picture in these figures shows a portion of $\mathbb{R}^2 \setminus D\mathbb{Z}^2$ near the origin. 

\begin{proof}
In the following argument, we assume, by Lemma \ref{Lem-change-spiral-numbers}, that $k_1 = 0$.

Now, suppose the two arcs are compatible.  Either $\sigma_1 = +$ or $\sigma_1 = -$.  
 
First, suppose $\sigma_1 = +$.  Then $C_1 = (m, 0, +)$.  Since $C_2$ is not isotopic to $C_1$ and lies entirely on one side of $C_1$, the only possibilities for $C_2$ are the ones shown in the left and middle pictures of Figure~\ref{Fig2a-2b-2f}. 
\begin{figure}
\begin{tabular}{ccc}
\scalebox{0.27}{\includegraphics{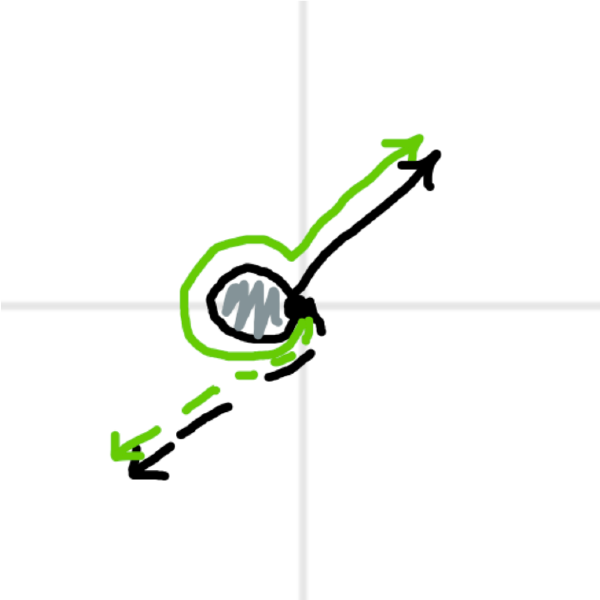}} &
\scalebox{0.27}{\includegraphics{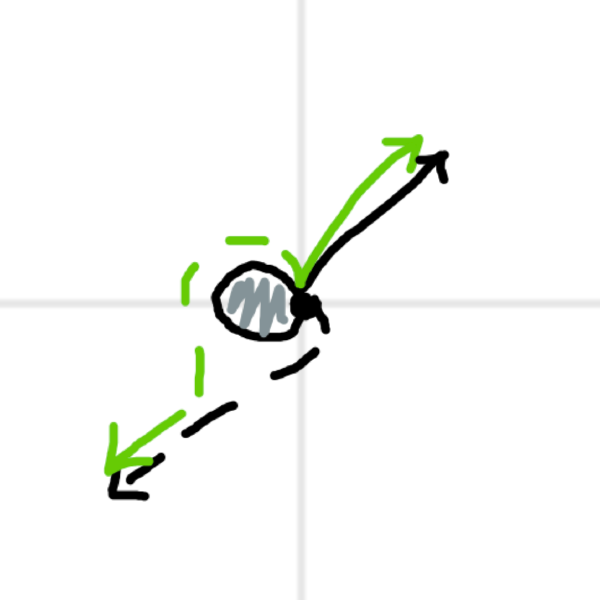}} &
\scalebox{0.27}{\includegraphics{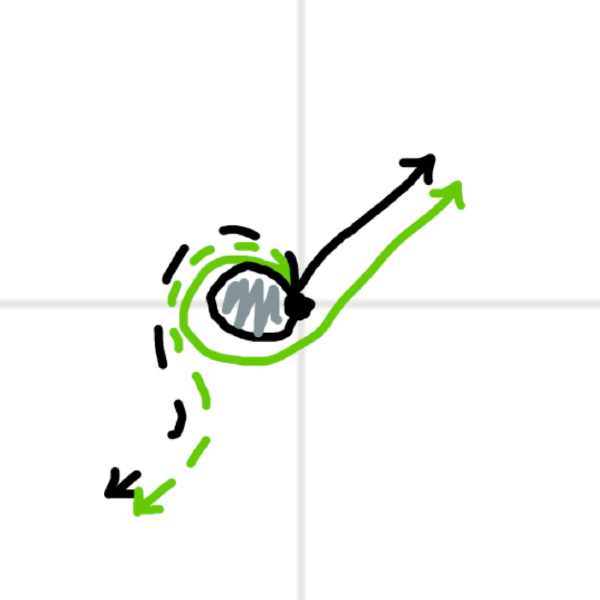}} \\
$C_1 = (1, 0, +)$ &
$C_1 = (1, 0, +)$ &
$C_1 = (1, 0, -)$ \\
$C_2 = (1, 1, -)$ &
$C_2 = (1, 0, -)$ &
$C_2 = (1, -1, +)$
\end{tabular}
\caption{Illustrations for the proof of Theorem \ref{Thm-compatible-arcs-same}. 
(Solid black line = $C_1$, solid green line = $C_2$, dashed black line = translate of~$C_1$ ending at origin, dashed green line = translate of $C_2$ ending at origin.).}
\label{Fig2a-2b-2f}
\end{figure} 
(Figure~\ref{Fig2a-2b-2f} shows the case where $m=1$, but changing the slope does not alter the argument.) 
In the left picture, $k_2=1$.  In the middle picture, $k_2 = 0$.  In both pictures, $\sigma_2 = -$.

Next, suppose $\sigma_1 = -$. If $k_2 = 0$, then $C_2 = (m,0,-) = C_1$ (a contradiction).  If $k_2 \neq 0$, the only possibility for $C_2$ is shown in the right picture of Figure~\ref{Fig2a-2b-2f} but it has $\sigma_2 = +$ (a contradiction).  

Conversely, suppose two arcs have standard lifts $C_1 = (m, 0, +)$ and $C_2 = (m, k_2, -)$ where $k_2 \in \{1, 0 \}$.
Both possibilities correspond to pairs of compatible arcs (see the left and middle pictures of Figure~\ref{Fig2a-2b-2f}).

The restatement of the condition in terms of the cover relation is immediate from the definition of the total order on arcs of slope $m$.
\end{proof}

To prove the second arc compatibility theorem, we make use of the following theorem.

\begin{thm}[\protect{\cite[Proposition 4.3]{28}}] \label{Prop-4-3-from-28}
Two arcs in the once-punctured torus whose standard lifts have non-origin endpoints $(a,b)$ and $(c,d)$ with $b/a < d/c$ are compatible if and only if $ad-bc=1$.
\end{thm}

\begin{thm} \label{Thm-compatible-arcs-distinct}
Two arcs in $(S,p)$ given by $C_1 = (b/a, k_1, \sigma_1)$ and $C_2 = (d/c, k_2, \sigma_2)$ with $b/a < d/c$ are compatible if and only if $ad - bc = 1$ and one of these holds:
\begin{enumerate}
    \item \label{Thm-compatible-arcs-distinct-part1} $\sigma_1 \leq \sigma_2$ and $k_2 = k_1$.
    \item \label{Thm-compatible-arcs-distinct-part2} $\sigma_1 > \sigma_2$ and $k_2 = k_1 + 1$.
\end{enumerate}
\end{thm}

\begin{proof}
In light of Theorem~\ref{Prop-4-1-from-28} and Theorem~\ref{Thm-arcs-as-images}, Theorem~\ref{Prop-4-3-from-28} tells us that the arcs given by $C_1$ and $C_2$ are compatible if and only if $ad-bc=1$ and their spirals around $D$ do not intersect (up to isotopy).  All that remains is to check that the given conditions on $\sigma_1$, $\sigma_2$, $k_1$, and $k_2$ precisely capture this non-intersection. 
We assume, by Lemma \ref{Lem-change-spiral-numbers}, that $k_1 = 0$.

First, suppose $\sigma_1 = +$.  As demonstrated in the left picture in Figure~\ref{Fig1a-1b}, $k_2$ is~$0$ or $1$.  
\begin{figure}
\begin{tabular}{cc}
\scalebox{0.28}{\includegraphics{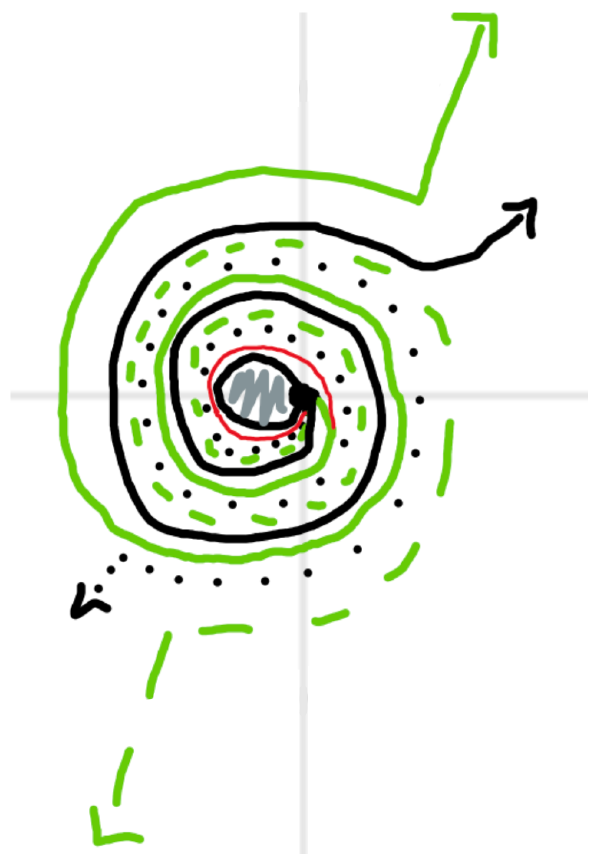}} &
\scalebox{0.28}{\includegraphics{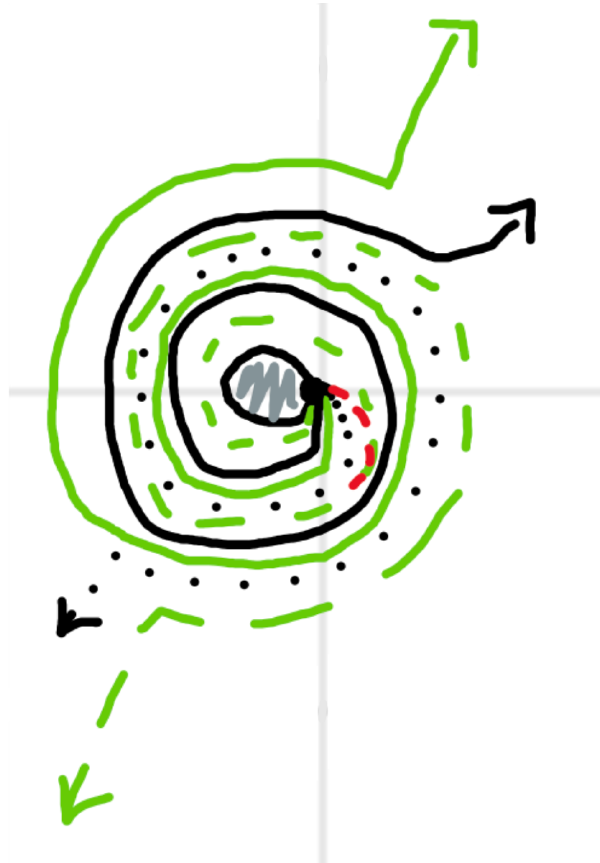}} \\
$C_1 = (1, 2, +)$ &
$C_1 = (1, 2, -)$ \\
$C_2$ is $(2, 2, +)$ or $(2, 3, -)$ &
$C_2 = (2, 2, \pm)$
\end{tabular}
\caption{Illustrations for the proof of Theorem \ref{Thm-compatible-arcs-distinct}. 
(Solid black line = $C_1$, solid green line = $C_2$, solid red line = alternate beginning for $C_2$, dashed black line = translate of $C_1$ ending at origin, dashed green line = translate of $C_2$ ending at origin, dashed red line = alternate end for $C_2$.).}
\label{Fig1a-1b}
\end{figure}
In both cases, the end of the arc with lift $C_2$ is forced to lie between the beginning and ending spirals of the arc with lift $C_1$.  However, when $k_2 = 0$, the beginning of the arc with lift $C_2$ lies CCW of the end of this arc (i.e., $\sigma_2 = +$) whereas, when $k_2 = 1$, the beginning of the arc with lift $C_2$ lies CW of the end of this arc (i.e., $\sigma_2 = -$).

On the other hand, suppose $\sigma_1 = -$.  As demonstrated in the right picture of Figure~\ref{Fig1a-1b}, $k_2=0$. 
However, as shown in the same picture, there are two ways in which the arc with lift $C_2$ can end.  Thus, $\sigma_2 = \pm$.  
\end{proof}

\begin{defn}[Cf. \protect{\cite[p. 89]{CATS1}}] \label{Defn-triangulation}
A \emph{triangulation} of $(S, p)$ or $(T,p)$ is a maximal collection of pairwise compatible arcs. 
A triangulation $\Delta$ divides $S$ or $T$ into topological triangles whose vertices are identified and whose sides are arcs in $\Delta$ and/or the boundary of $S$. 
(In some marked surfaces, the same arc can serve as multiple sides of a given triangle. However, this does not happen in $(S,p)$ or $(T,p)$; see Theorem~\ref{thm-triangles-only-slope} and \cite[Section~3]{28}.) 
\end{defn}

The following theorem characterizes the triangulations of $(S, p)$. 

\begin{thm} \label{thm-triangulations}
Each triangulation of $(S,p)$ consists of four arcs. 
The arcs given by $(\frac{b_1}{a_1}, k_1, \sigma_1)$, $(\frac{b_2}{a_2}, k_2, \sigma_2)$, $(\frac{b_3}{a_3}, k_3, \sigma_3)$, and $(\frac{b_4}{a_4}, k_4, \sigma_4)$ with $\frac{b_1}{a_1} \leq \frac{b_2}{a_2} \leq \frac{b_3}{a_3} \leq \frac{b_4}{a_4}$ form a triangulation if and only if exactly one equality holds (say, $\frac{b_n}{a_n} = \frac{b_{n+1}}{a_{n+1}}$), $b_i a_j - a_i b_j = 1$ for all $i < j$ such that $i, j \in \{ 1,2,3,4 \}$ with $(i, j) \neq (n, n+1)$, and either
\begin{enumerate}
     \item $\sigma_1 \leq \sigma_2 \leq \sigma_3 \leq \sigma_4$ with exactly one strict inequality holding (namely, $\sigma_n < \sigma_{n+1}$) and $k_1 = k_2 = k_3 = k_4$ or
     \item $\sigma_1 \geq \sigma_2 \geq \sigma_3 \geq \sigma_4$ with exactly one strict inequality holding (namely, $\sigma_n > \sigma_{n+1}$) and  $k_1 \leq k_2 \leq k_3 \leq k_4$ with exactly one strict inequality holding (namely, $k_{n+1} = k_n + 1$).
    \end{enumerate}
\end{thm}

\begin{proof}
A triangulation is a maximal collection of pairwise compatible arcs. All triangulations of a given surface have the same number of arcs in them \cite[Lemma~3.6]{CATS1}.  To see that this number is $4$ for $(S,p)$, observe that the arcs given by $(-1,0,-)$, $(-1, 0, +)$, $(0, 0, +)$, and $(\infty, 0, +)$ divide $S$ into triangles and therefore form a triangulation.  See Figure~\ref{FigforTriangulations}.
\begin{figure}
\scalebox{0.38}{\includegraphics{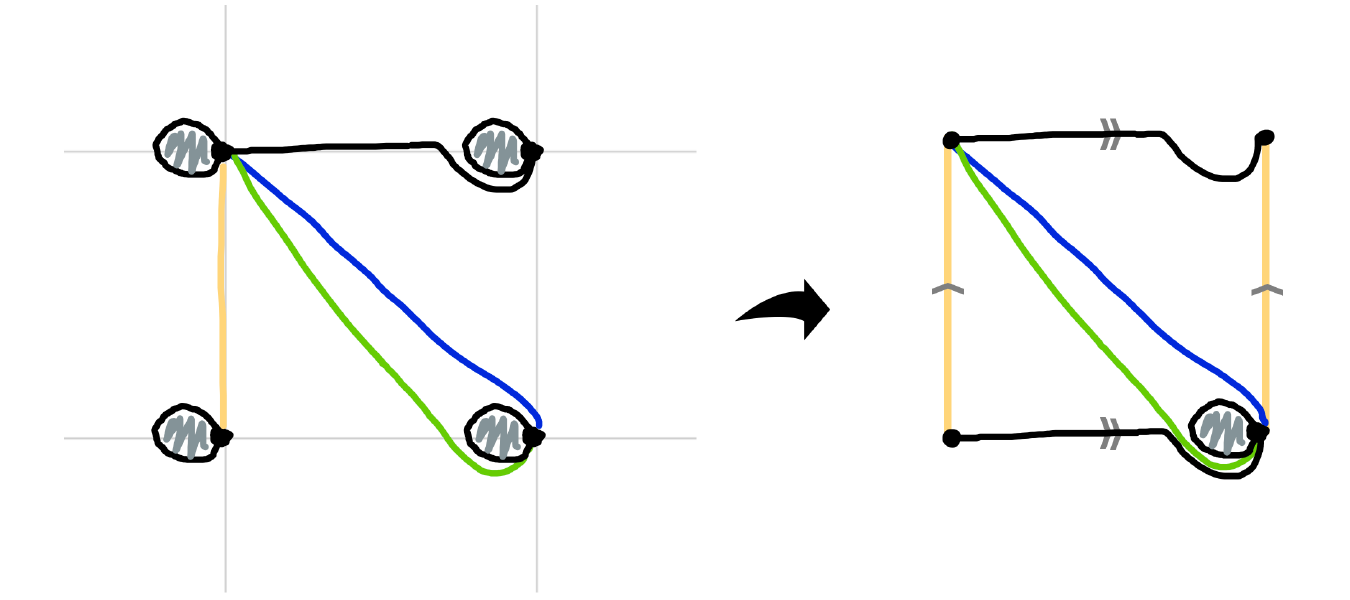}}
\caption{\emph{Left:}  In $\mathbb{R}^2 \setminus D\mathbb{Z}^2$, the standard lifts $(-1,0,-)$ (blue), $(-1, 0, +)$ (green), and $(0, 0, +)$ (black) as well as the lift of the arc given by $(\infty, 0, +)$ that ends at the origin (yellow).  \emph{Right:}  The corresponding arcs in $(S,p)$.}
\label{FigforTriangulations}
\end{figure}

Now, consider four arcs in $(S,p)$ that form a triangulation.  
They have at most three distinct slopes. 
(This is easily seen by comparing Theorem~\ref{Prop-4-3-from-28} and Theorem~\ref{Thm-compatible-arcs-distinct} and noting that  triangulations of the once-punctured torus consist of three arcs.) 
Theorem~\ref{Thm-compatible-arcs-same} implies that a given slope can appear at most twice in a triangulation.  Thus, either exactly one pair of arcs have the same slope or exactly two pairs of arcs have the same slope.

To see that the latter is impossible, let $C_1 = (m_1, k_1, \sigma_1)$, $C_1'=(m_1,k_1',\sigma_1')$, $C_2=(m_2, k_2, \sigma_2)$, and $C_2'=(m_2,k_2', \sigma_2')$ with $\sigma_1 \geq \sigma_1'$, $\sigma_2 \geq \sigma_2'$, and $m_1 < m_2$ be the standard lifts of arcs forming a triangulation.  
Theorem~\ref{Thm-compatible-arcs-same} implies that $\sigma_1 = \sigma_2 = +$ and $\sigma_1' = \sigma_2' = -$.  
By Theorem~\ref{Thm-compatible-arcs-distinct}, since the arcs given by $C_1$ and~$C_2$ are compatible and the arcs given by $C_1'$ and $C_2$ are compatible, $k_2 = k_1$ and $k_2 = k_1'$.  However, since the arcs given by $C_1$ and $C_2'$ are compatible and the arcs given by $C_1'$ and $C_2'$ are compatible, $k_2' = k_1 + 1$ and $k_2' = k_1'$ (a contradiction).  Thus, any triangulation has exactly one pair of arcs with the same slope.  

Next, consider a triangulation of $(S,p)$ whose arcs are given by $C_1 = (\frac{b_1}{a_1}, k_1, \sigma_1)$, $C_2 = (\frac{b_2}{a_2}, k_2, \sigma_2)$, $C_3 = (\frac{b_3}{a_3}, k_3, \sigma_3)$, and $C_4 = (\frac{b_4}{a_4}, k_4, \sigma_4)$ with ${\frac{b_1}{a_1} \leq \frac{b_2}{a_2} \leq \frac{b_3}{a_3} \leq \frac{b_4}{a_4}}$ and with exactly one 
equality holding (say, $\frac{b_n}{a_n} = \frac{b_{n+1}}{a_{n+1}}$). 
Theorem~\ref{Thm-compatible-arcs-distinct} says that $b_i a_j - a_i b_j = 1$ for all $i < j$ where $i, j \in \{ 1,2,3,4 \}$ and $(i, j) \neq (n, n+1)$. 
By Theorem~\ref{Thm-compatible-arcs-same}, $C_n$ and $C_{n+1}$ have opposite sign and we can renumber, if necessary, so that $\sigma_n < \sigma_{n+1}$. 
By Theorem~\ref{Thm-compatible-arcs-same}, either $k_n = k_{n+1}$ or $k_{n+1} = k_n+1$. 

Suppose $k_n = k_{n+1}$. 
We use Theorem~\ref{Thm-compatible-arcs-distinct} repeatedly. 
Since the arcs given by $C_n$ and $C_i$ for all $i > n+1$ are compatible and $\sigma_n = -$, $k_n = k_i$ for all $i > n+1$.  However, since $k_n = k_{n+1}$, this means that $k_{n+1} = k_i$ for all $i > n+1$.  Since the arcs given by~$C_{n+1}$ and $C_i$ for all $i > n+1$ are compatible, $\sigma_{n+1} \leq \sigma_i$ for all $i > n+1$.  
Similarly, since the arcs given by $C_j$ and $C_{n+1}$ for all $j < n$ are compatible and $\sigma_{n+1} = +$, $k_j = k_{n+1}$ for all $j < n$.  However, since $k_n = k_{n+1}$, this means that $k_j = k_n$ for all $j < n$.  Since the arcs given by $C_j$ and $C_n$ for all $j < n$ are compatible, $\sigma_j \leq \sigma_n$ for all $j < n$.  Thus, we see that $\sigma_1 \leq \sigma_2 \leq \sigma_3 \leq \sigma_4$ and $k_1 = k_2 = k_3 = k_4$. 

Suppose $k_{n+1} = k_n+1$. 
We again use Theorem~\ref{Thm-compatible-arcs-distinct} repeatedly.  
Since the arcs given by $C_{n+1}$ and $C_i$ for all $i > n+1$ are compatible and $\sigma_{n+1} = -$, $k_{n+1}=k_i$ for all $i > n+1$.
However, since $k_{n+1} = k_n+1$, this means that $k_i = k_n+1$ for all $i > n+1$.  
Since the arcs given by $C_n$ and $C_i$ for all $i > n+1$ are compatible, $\sigma_n > \sigma_i$ for all $i > n+1$.  
Similarly, since the arcs given by $C_j$ and $C_n$ for all $j < n$ are compatible and $\sigma_n = +$, $k_j = k_n$ for all $j < n$.  
 However, since $k_{n+1} = k_n+1$, this means that $k_{n+1} = k_j + 1$ for all $j < n$.  Since the arcs given by $C_j$ and~$C_{n+1}$ for all $j < n$ are compatible, $\sigma_j > \sigma_{n+1}$ for all $j < n$.  Thus, we see that $\sigma_1 \geq \sigma_2 \geq \sigma_3 \geq \sigma_4$ and $k_1 \leq k_2 \leq k_3 \leq k_4$ (with equality except that $k_{n+1} = k_n+1$). 

Conversely, suppose there are four arcs with standard lifts $C_1 = (\frac{b_1}{a_1}, k_1, \sigma_1)$, $C_2 = (\frac{b_2}{a_2}, k_2, \sigma_2)$, $C_3 = (\frac{b_3}{a_3}, k_3, \sigma_3)$, and $C_4 = (\frac{b_4}{a_4}, k_4, \sigma_4)$ with $\frac{b_1}{a_1} \leq \frac{b_2}{a_2} \leq \frac{b_3}{a_3} \leq \frac{b_4}{a_4}$ and exactly one equality holding (say, $\frac{b_n}{a_n} = \frac{b_{n+1}}{a_{n+1}}$), $b_i a_j - a_i b_j = 1$ for all $i < j$ where $i, j \in \{ 1,2,3,4 \}$ with $(i, j) \neq (n, n+1)$, and either:
\begin{enumerate}
     \item $\sigma_1 \leq \sigma_2 \leq \sigma_3 \leq \sigma_4$ with exactly one strict inequality holding (namely, $\sigma_n < \sigma_{n+1}$) and $k_1 = k_2 = k_3 = k_4$ or
     \item $\sigma_1 \geq \sigma_2 \geq \sigma_3 \geq \sigma_4$ with exactly one strict inequality holding (namely, $\sigma_n > \sigma_{n+1}$) and  $k_1 \leq k_2 \leq k_3 \leq k_4$ with exactly one strict inequality holding (namely, $k_{n+1} = k_n + 1$).
    \end{enumerate}
Direct applications of Theorems \ref{Thm-compatible-arcs-same} and  \ref{Thm-compatible-arcs-distinct} show that any two of these four arcs are compatible with each other.
\end{proof}

It is immediate from Definition~\ref{Defn-triangulation} that each triangle resulting from a triangulation of $(S,p)$ is bounded by either (a) three arcs or (b) two arcs and the boundary of $S$.  The following theorem characterizes these triangles.

\begin{thm} \label{thm-triangles-only-slope}
Let $\Delta$ be a triangulation of $(S,p)$. 
Three arcs in $\Delta$ bound a triangle if and only if their slopes are distinct. Two arcs in $\Delta$, along with the boundary of~$S$, bound a triangle if and only if the arcs' slopes are equal. 
The triangulation $\Delta$ divides~$S$ into three triangles.
\end{thm}

\begin{proof}
Suppose three arcs in $\Delta$ bound a triangle.  
If two of these arcs have the same slope, then these two arcs divide $S$ into two connected regions: (a) a region containing the boundary and (b) an annular region whose two boundary circles are the two arcs (and so share a point, namely, $p$).  The third arc cannot lie in the annular region since, in the annulus, it would be isotopic to each of the other two arcs.  However, it cannot lie in the other region, for then it would be isotopic to one of the other arcs.  Therefore, the slopes are distinct. 
 
Conversely, suppose three arcs in $\Delta$ have distinct slopes.  
To see that these arcs bound a connected region, observe that any arc cuts $S$ into an annulus with the boundary of $S$ touching one boundary circle of the annulus.  
A second arc, compatible and with a different slope, cuts the annulus into a quadrilateral with the boundary of $S$ lying at one corner.  
A third arc, compatible with the other two arcs and with a different slope from either of them, divides this quadrilateral into two triangles, one of which has the boundary of $S$ at a corner.  The triangle not containing the boundary of $S$ is bounded by the three arcs.

Next, suppose two arcs in $\Delta$, along with the boundary of $S$, bound a triangle.  Then one of the arcs is isotopic to a curve that begins at $p$, travels closely around the boundary of $S$, and then follows along the other arc until it returns to $p$.  Thus, both arcs have the same slope.

Conversely, suppose two arcs in $\Delta$ have the same slope.  By Theorem~\ref{Thm-compatible-arcs-same}, these arcs are given by $(m, k_1, +)$ and $(m, k_2, -)$ where $k_2 \in \{k_1, k_1+1 \}$.  
In either case ($k_2=k_1$ or $k_2=k_1+1$), the two arcs, along with the boundary of $S$, bound a triangle (as demonstrated in the left and middle pictures of Figure~\ref{Fig2a-2b-2f}). 

From this argument, we see that any subset of $\Delta$ consisting of three arcs with distinct slopes or two arcs with the same slope corresponds to a unique triangle. 
By Theorem~\ref{thm-triangulations}, $\Delta$ consists of four arcs\textemdash exactly two of which have the same slope. 
Thus, there are exactly two subsets of $\Delta$ consisting of three arcs with distinct slopes and exactly one subset of $\Delta$ consisting of two arcs of the same slope. 
Hence, $\Delta$ divides $S$ into $2+1=3$ triangles. 
\end{proof}

\section{Allowable Curves}

In this section, we characterize the dread torus' allowable curves, its maximal collections of compatible allowable curves, and the shear coordinates of its allowable curves. 
These characterizations will be used later in proving that the dread torus has the Null Tangle property and describing the dread torus' mutation fan. 

\begin{defn}[Cf. \protect{\cite[Definitions 3.7 and 4.1]{27}}] \label{Defn-allowable-curve}
In $(S,p)$, an \emph{allowable curve} is a non-self-intersecting curve that is non-contractible and is either (a) a closed curve (i.e., embedded circle) lying in the interior of $S$ or (b) a curve whose two endpoints are distinct, unmarked boundary points but which otherwise lies in the interior of~$S$ and is not contractible to part of the boundary.  
The closed allowable curves in $(T, p)$ are defined the same way except they cannot intersect $p$ or contract to~$p$. 
(We will not need the definition of non-closed allowable curves in $(T,p)$.) 
Allowable curves are considered up to isotopy fixing $p$. 
\end{defn}

To characterize the allowable curves of $(S,p)$, we make use of a map $\kappa$ defined in~\cite{27}.    
Let $\alpha$ be an arc.  Then, in the dread torus, $\kappa(\alpha)$ is defined, up to isotopy, as the curve that is the same as~$\alpha$ except that each endpoint is an unmarked boundary point slightly clockwise of $p$ (i.e., as the curve approaches $p$, it misses $p$ slightly to the left).
The two endpoints of $\kappa (\alpha)$ are slightly different distances from $p$ with the distances chosen so that $\kappa(\alpha)$ does not self-intersect.  
Lemma 5.1 in \cite{27}, as specialized to $(S,p)$, states that $\kappa$ is a bijection between the arcs of $(S,p)$ and its non-closed allowable curves.  Thus, since the arcs of $(S,p)$ are in bijection with $\{ (m, k, \sigma) : m \in \mathbb{Q} \cup \{ \infty \}, k \in \mathbb{Z}, \sigma \in \{ +, - \} \}$, there is a bijection between these triples and the non-closed allowable curves of $(S,p)$.  If an arc $\alpha$ is given by $(m, k, \sigma)$, we say that~$\kappa(\alpha)$ is given by $(m, k, \sigma)$ and has slope $m$, $k$ spirals, and sign~$\sigma$. 
We use the same total order on non-closed allowable curves having the same slope as we did on arcs having the same slope.

\begin{thm} \label{ACs}
The allowable curves in $(S,p)$ can be classified as follows.
\begin{description}
\item[\textbf{Type~1 curves}]  One non-closed allowable curve for each triple $(m, k, \sigma)$, obtainable by applying the map $\kappa$ to the arc given by $(m, k, \sigma)$.  
\item[\textbf{Type~2 curves}]  One closed allowable curve for each $m \in \mathbb{Q} \cup \{ \infty \}$, obtainable by shifting both endpoints of any arc with slope $m$ to an interior point of $S$ near $p$.
\item[\textbf{Type~3 curve}]  The closed allowable curve isotopic to the boundary of $S$.
\end{description}
\end{thm}

We say that a Type~2 curve obtained from an arc of slope $m$ has slope $m$.

\begin{proof}
The existence of a bijection between non-closed allowable curves and triples and the means of obtaining a non-closed allowable curve from an arc were justified earlier.  

Now, we show there is exactly one closed allowable curve for each $m \in \mathbb{Q} \cup \{ \infty \}$.
For each triple $(m, k, \sigma)$, define $\lambda(m, k, \sigma)$ to be the curve that is formed from the arc given by $(m, k, \sigma)$ by shifting both of its endpoints to a point $x$ near~$p$ in the interior of $S$.  By construction, $\lambda(m, k, \sigma)$ is a closed allowable curve.

To show that $\lambda(m, k, \sigma)$ is the same allowable curve regardless of the values of~$k$ and $\sigma$, we first show that $\lambda(m, k, \sigma)$ for $k \neq 0$ is isotopic to $\lambda(m, 0, \sigma)$ and then show that $\lambda(m, 0, +)$ and $\lambda(m, 0, -)$ are isotopic.  
Let $C$ be the standard lift of the arc given by $(m, k, \sigma)$ ($k \neq 0$), and let $C'$ be the lift of that arc that ends at the origin.   
As illustrated in Figure \ref{k-and-sign}, the curves $C$ and $C'$ stay close together as they spiral out from the origin 
and, up to isotopy, only differ in that one of them completes an additional, final half spiral. 
\begin{figure}
\begin{tabular}{cccc}
\scalebox{.8}{\includegraphics{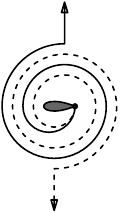}} &
\scalebox{.8}{\includegraphics{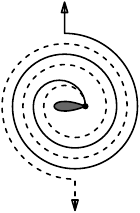}} &
\scalebox{.8}{\includegraphics{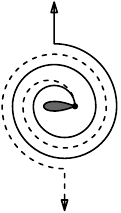}} &
\scalebox{.8}{\includegraphics{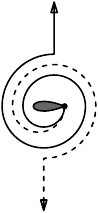}} \\
$(\infty, 2, +)$  &
$(\infty, -2, -)$ &
$(\infty, -2, +)$  &
$(\infty, 2, -)$ \\
\end{tabular}
\caption{Illustrations for the proof of Theorem~\ref{ACs} (solid black line = $C$; dashed black line = $C'$).}
\label{k-and-sign}
\end{figure} (If $k$ and $\sigma$ have the same parity, 
$C'$ completes the additional half spiral; if $k$ and $\sigma$ have opposite parity, 
$C$ completes the additional half spiral.)   
Thus, there is an isotopy of $\lambda(m,k,\sigma)$ that removes all its spirals except for a half spiral.  (Starting from the beginning of the non-spiraling portion of $\lambda(m,k,\sigma)$, the half spiral travels CW around the boundary of $S$ if $\sigma = +$ and travels CCW around the boundary of $S$ if $\sigma = -$.)  Hence, $\lambda(m,k,\sigma)$ is isotopic to $\lambda(m,0,\sigma)$. 
Notice that $\lambda(m, 0, +)$ and $\lambda(m, 0, -)$ divide $S$ into two connected regions:  (a) a region containing the boundary of $S$ and (b) an annular region whose two boundary circles are $\lambda(m, 0, +)$ and $\lambda(m, 0, -)$.  Since there is an isotopy between the two boundary circles of an annulus, $\lambda(m, 0, +)$ and $\lambda(m, 0, -)$ are isotopic.  

To see that $\lambda(b/a, k_1, \sigma_1)$ and $\lambda(d/c, k_2, \sigma_2)$ are not isotopic if $b/a \neq d/c$, notice that $\lambda(b/a, k_1, \sigma_1)$ wraps around the longitude of $S$ a total of $a$ times and around the latitude of $S$ a total of $b$ times.  
Since $d/c \neq b/a$, $d \neq b$ or $c \neq a$.  Thus, $\lambda(d/c, k_2, \sigma_2)$ wraps around either the longitude or the latitude of $S$ a different number of times than $\lambda(b/a, k_1, \sigma_1)$ does.  Therefore, $\lambda(b/a, k_1, \sigma_1)$ and $\lambda(d/c, k_2, \sigma_2)$ intersect.

Next, observe that a circle isotopic to the boundary of $S$ is a closed allowable curve.  

Finally, we show that a closed curve $C$ that is not isotopic to the boundary of~$S$ is a Type~2 curve.  
Let~$\gamma$ denote a non-self-intersecting path in $(S,p)$ from $p$ to~$C$ that is disjoint from the boundary of $S$ and $C$ except at its endpoints.  Let $\alpha$ be a path that begins at $p$, travels just right of $\gamma$, turns right onto $C$, traverses $C$ once, and then turns right and returns to $p$ by traveling along $\gamma$.  Then $\alpha$ lies in the interior of $S$ except for its endpoint $p$, is non-self-intersecting, is non-contractible (else, $C$ would be contractible), and is not isotopic to the boundary of $S$ (else, $C$ would be isotopic to the boundary of $S$).  Thus, $\alpha$ is an arc with a standard lift $(m, k, \sigma)$.  By the construction of $\alpha$, $\lambda(m, k, \sigma)$ is isotopic to $C$.  Therefore, $C$ is a Type~2 curve.
\end{proof}

\begin{rem} \label{f_n-applied-to-ACs}
The homeomorphism $f_n$ of $S$ described in the proof of Lemma~\ref{Lem-change-spiral-numbers} increases the number of spirals of every Type~1 curve by $n$ and does not, up to isotopy, otherwise alter any of the allowable curves. 
\end{rem}

\begin{defn}[Cf. \protect{\cite[Definition 4.1]{27}}] 
Two allowable curves in $(S,p)$ are \emph{compatible} if they have isotopy representatives that do not intersect. 
\end{defn}

To classify the maximal collections of pairwise compatible allowable curves, we make use of the following lemma.

\begin{lem}[\protect{\cite[Part 2 of Lemma 5.1]{27}}] \label{art-27-lem-5.1}
In any marked surface, two arcs $\alpha$ and~$\beta$ are compatible if and only if the allowable curves $\kappa (\alpha)$ and $\kappa (\beta)$ are compatible. 
\end{lem}

\begin{thm} \label{Thm-maxlcolsofACs}
In $(S,p)$, the maximal collections of pairwise compatible allowable curves are as follows.
\begin{description}
\item[Type~A (Type~1 curves only)]  The triangulations of $S$ with each arc $\alpha$ replaced by the non-closed allowable curve $\kappa(\alpha)$. 
\item[Type~B (Type~1 \& 2 curves mixed)]  All sets consisting of a closed allowable curve of slope $m$ and two non-closed allowable curves given by $(m, k_1, +)$ and $(m, k_2, -)$ where $k_1$ is any fixed integer and $k_2 \in \{ k_1, k_1+1 \}$.
\item[Type~C (Type~2 \& 3 curves mixed)]  All sets consisting of a closed allowable curve of slope $m$ and the Type~3 curve. 
\end{description} 
\end{thm}

\begin{figure}
\scalebox{0.29}{\includegraphics{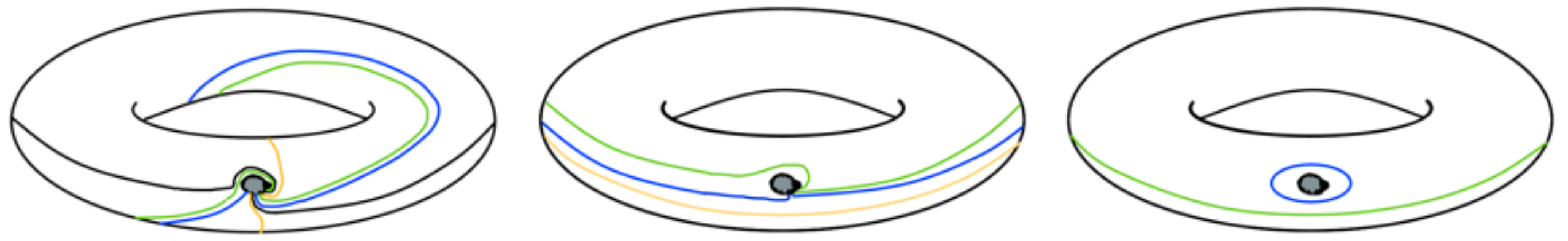}}
\caption{The three types of maximal collections of allowable curves.  From left to right: Type A, Type B, and Type C.}
\label{maxl-cols-of-ACs}
\end{figure}

This theorem is illustrated in Figure \ref{maxl-cols-of-ACs}.  

\begin{proof}
Two facts are immediate: any Type~1 curve is incompatible with the Type~3 curve and any Type~2 curve is compatible with the Type~3 curve.  
In addition, any two Type~2 curves are incompatible and a Type~1 curve is compatible with a Type~2 curve if and only if they have the same slope 
(if two allowable curves have different slopes, they wrap around either the longitude or the latitude of $S$ a different number of times; if a Type~1 curve and a Type~2 curve have the same slope, fix a representative for the Type~1 curve whose spirals are close to the boundary of~$S$ and then choose a representative for the Type~2 curve that runs closely along the Type~1 curve's non-spiraling portion and circumvents its spirals). 
The theorem follows from these four facts, the definition of a triangulation, Lemma~\ref{art-27-lem-5.1}, and Theorem~\ref{Thm-compatible-arcs-same}. 
\end{proof}

\begin{rem} \label{Rem-cardinality-of-maxl-cols}
Each Type~A collection consists of $4$ curves, each Type~B collection consists of $3$ curves, and each Type~C collection consists of $2$ curves.   
\end{rem}

\begin{defn}[Cf. \protect{\cite[Definition 12.2]{CATS2}}] 
Let $\Delta = \{ \alpha_1, \alpha_2, \alpha_3, \alpha_4 \}$ be a triangulation of $(S,p)$, and let $\lambda$ be an allowable curve. 
The \emph{shear coordinates} (or \emph{shear coordinate  vector}) of $\lambda$ with respect to $\Delta$, denoted $\mathbf{b}(\Delta, \lambda)$, is indexed by the arcs of $\Delta$. 
Each entry $b_{\alpha_i}(\Delta, \lambda)$ is a sum whose terms come from $\{ 1, -1, 0 \}$.  
Each term corresponds to a crossing of $\alpha_i$ by $\lambda$, and the value of this term is determined by how $\lambda$ intersects the quadrilateral formed by the two triangles that share $\alpha_i$ as a side.  
Every time $\lambda$ enters and exits this quadrilateral through opposite sides, the term is $\pm1$. 
It is -1 if the opposite sides form the top and bottom of a ``Z'' shape that has $\alpha_i$ as its diagonal; it is 1 if the opposite sides form the top and bottom of a ``backwards Z'' shape that has $\alpha_i$ as its diagonal.  
Every time $\lambda$ enters and exits the quadrilateral through adjacent sides, the term is $0$.  
Up to isotopy, $\lambda$ never crosses $\alpha_i$ in one direction and then in the opposite direction, so $b_{\alpha_i}(\Delta, \lambda)$ is either a sum of nonnegative numbers or a sum of nonpositive numbers.
\end{defn}

Let $\Delta_0$ denote the triangulation of $(S,p)$ consisting of the arcs $\gamma_1$, $\gamma_2$, $\gamma_3$, and~$\gamma_4$ given by $(0, 0, +)$, $(\infty, 0, +)$, $(-1, 0, +)$, and $(-1, 0, -)$, respectively (shown in Figure~\ref{FigforTriangulations}).  
We calculate the shear coordinates of the dread torus' allowable curves with respect to $\Delta_0$. 
To make these calculations, we use the following lemma. 

\begin{lem} \label{lifts-of-ACs}
In $(S,p)$: 
\begin{enumerate}
\item \label{lifts-of-ACs-T1} The Type~1 curve given by $(b/a, k, \sigma)$ lifts to a curve in $\mathbb{R}^2 \setminus D\mathbb{Z}^2$ identical to the standard lift of the arc given by $(b/a, k, \sigma)$ except that its endpoints are different (they are boundary points slightly clockwise of $(0,0)$ and $(a,b)$). 

\item \label{lifts-of-ACs-T2}  The Type~2 curve with slope $b/a$ lifts to a line in $\mathbb{R}^2 \setminus D\mathbb{Z}^2$ with slope $b/a$ that passes through a fixed point $P$ nearby $(0,0)$ that lies between $(0,0)$ and $(1,0)$ if $b/a \neq 0$ and lies between $(0,0)$ and $(0,-1)$ if $b/a = 0$.
 
\item \label{lifts-of-ACs-T3} The Type~3 curve lifts to a circle in $\mathbb{R}^2 \setminus D\mathbb{Z}^2$ around the boundary component at $(0,0)$.
\end{enumerate}
\end{lem}

The lifted curve or line in parts (\ref{lifts-of-ACs-T2}) or (\ref{lifts-of-ACs-T3}) of this lemma is only defined if the size of the removed disks of $\mathbb{R}^2\setminus D\mathbb{Z}^2$ is sufficiently small. 
Thus, in applications of this lemma, we either fix a size for the removed disks that is ``small enough'' for all the Type~1 and/or~2 curves under consideration (if the curves correspond to a finite number of slopes) or permit the lifts to bend slightly around boundary components of $\mathbb{R}^2 \setminus D\mathbb{Z}^2$ (if the curves correspond to an infinitude of slopes). 
Viewing the boundary of $S$ as a ``thickened point'' is, in essence, fixing a size for the removed disks that is sufficiently small for \emph{all} the Type~1 and~2 curves of $(S,p)$. 

For each allowable curve in $(S,p)$, we refer to its lift described in Lemma~\ref{lifts-of-ACs} as its \emph{standard lift}. 

\begin{proof}
Statements (\ref{lifts-of-ACs-T1}) and (\ref{lifts-of-ACs-T3}) follow directly from the constructions of the Type~1 curves and the Type~3 curve described in Theorem~\ref{ACs}. 
For Statement~(\ref{lifts-of-ACs-T2}), it is clear from 
paragraphs 2 and 3 
of the proof of Theorem~\ref{ACs} that the desired lift can be obtained by performing the following operations to the standard lift of the arc given by $(b/a, 0, \sigma)$ where $\sigma = +$ if $b/a > 0$ and $\sigma = -$ if $b/a \leq 0$: first, translate the lift in the direction perpendicular to the lift's straight part and toward the lift's partial spirals (by an amount sufficient for the following operation to be performed); next, ``straighten out'' the lift's partial spirals (so the lift becomes a line segment of slope~$b/a$); and, finally, extend the lift to be a full line. (The last is necessary since an arc is an embedded closed interval but a Type~2 curve is an embedded circle.) 
\end{proof}
\begin{thm} \phantomsection\label{Thm-shear-coords-of-ACs}
\begin{enumerate} 
\item The shear coordinates for the Type~3 curve are $[0, 0, 1, -1]$.

\item The shear coordinates for a Type~2 curve given by $b/a$ are: 

    $[-b, a, b, -a]$ if $ 0 \leq b/a \leq \infty $,
    
    $[-b, a+2b, 0, -a-b]$ if $ -1 \leq b/a \leq 0$, or
    
    $[b+2a, -a, -b-a, 0]$ if $ - \infty < b/a \leq -1 $.

\item \label{Thm-shear-coords-ACs-slope-greater-than-0} The shear coordinates for a Type~1 curve given by $(b/a, k, \sigma)$ with \break $0 < b/a < \infty$ are:

    $[-b-1, a-1, b-2k+1, -a+2k]$ if $k < 1$ and $\sigma = -$.

    $[-b-1, a-1, b-2k, -a+2k+1]$ if $k < 1$ and $\sigma = +$.

    $[-b+1,a+1,b+2k-4,-a-2k+3]$ if $k > 1$ and $\sigma = -$. 

    $[-b+1,a+1,b+2k-3,-a-2k+2]$ if $k \geq 1$ and $\sigma = +$. 
    
    $[-b, a, b-1, -a+1]$ if $k = 1$ and $\sigma = -$.

\item \label{Thm-shear-coords-ACs-slope-between-0-and-neg-1}The shear coordinates for a Type~1 curve given by ${(b/a, k, \sigma)}$ with \break $-1 < b/a < 0$ are:

    $[-b-1, a+2b+1, -2k, -a-b+2k-1]$ if $k < 1$ and $\sigma = -$. 

    $[-b-1, a+2b+1, -2k-1, -a-b+2k]$ if $k < 1$ and $\sigma = +$.

    $[-b+1,a+2b-1,2k-3,-a-b-2k+4]$ if $k > 1$ and $\sigma = -$.

    $[-b+1,a+2b-1,2k-2,-a-b-2k+3]$ if $k \geq 1$ and $\sigma = +$.
    
    $[-b, a+2b, -1, -a-b+1]$ if $k = 1$ and $\sigma = -$.
    
\item \label{Thm-shear-coords-ACs-slope-less-than-neg-1} The shear coordinates for a Type~1 curve given by $(b/a, k, \sigma)$ with \break $- \infty < b/a < -1$ are:

    $[2a+b, -a, -a-b-1, 1]$ if $k = 0$ and $\sigma = +$.

   $[2a+b+1, -a-1, -a-b-2k-1, 2k]$ if $k \leq 0$ and $\sigma = -$.

   $[2a+b+1, -a-1, -a-b-2k-2, 2k+1]$ if $k < 0$ and $\sigma = +$.

   $[2a+b-1, -a+1, -a-b+2k-2, -2k+3]$ if $k > 0$ and $\sigma = -$.

   $[2a+b-1, -a+1, -a-b+2k-1, -2k+2]$ if $k > 0$ and $\sigma = +$.

\item \label{Thm-shear-coords-ACs-slope-equals-0} The shear coordinates for a Type~1 curve given by $(0, k, \sigma)$ can be obtained by setting $b=0$ and $a=1$ in  Case~\ref{Thm-shear-coords-ACs-slope-greater-than-0} if $k<1$ and in Case~\ref{Thm-shear-coords-ACs-slope-between-0-and-neg-1} if $k \geq 1$. 

\item \label{Thm-shear-coords-ACs-slope-equals-neg-1} The shear coordinates for a Type~1 curve given by $(-1, k, \sigma)$ can be obtained by setting $b=-1$ and $a=1$ in Case~\ref{Thm-shear-coords-ACs-slope-between-0-and-neg-1} if $k<1$ and in Case~\ref{Thm-shear-coords-ACs-slope-less-than-neg-1} if $k \geq 1$.

\item \label{Thm-shear-coords-ACs-slope-equals-infty} The shear coordinates for a Type~1 curve given by $(\infty, k, \sigma)$ can be obtained by setting $b=1$ and $a=0$ in Case~\ref{Thm-shear-coords-ACs-slope-greater-than-0} if $k \geq 1$; setting $b=-1$, $a=0$, and $\sigma=-\sigma$ in Case~\ref{Thm-shear-coords-ACs-slope-less-than-neg-1} if $k < 1$ and $\sigma=+$; and setting $b=-1$, $a=0$, $\sigma=-\sigma$, and $k=k-1$ in Case~\ref{Thm-shear-coords-ACs-slope-less-than-neg-1} if $k < 1$ and $\sigma=-$.

\end{enumerate}
\end{thm}

The full proof of Theorem~\ref{Thm-shear-coords-of-ACs} consists mostly of simple, but tedious, calculations.  The details are omitted here, due to length, but can be found in \cite{mydissertation}. 

\begin{proof}[Proof sketch.]
\emph{Case 1.} The Type~3 curve's shear coordinates can be calculated directly from Figure~\ref{Type-3-AC-shear-coords}. 
\begin{figure}
\scalebox{.8}{\includegraphics{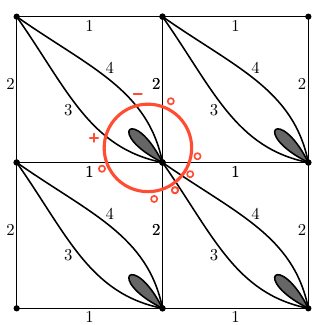}}
\caption{The Type~3 curve's standard lift and the lifts of $\Delta_0$ that intersect it.}
\label{Type-3-AC-shear-coords}
\end{figure}

\emph{Case 2.}  First, suppose $b/a > 0$. 
If we view the boundary of $S$ as a ``thickened point,'' the standard lift of the Type~2 curve of slope $b/a$ becomes the lift $l$ of a closed allowable curve $\lambda$ in $(T,p)$ by Definition~\ref{Defn-allowable-curve}. 
In \cite{28}, it is shown that the arc crossings of $\lambda$ can be encoded in a word formed by recording a letter ($\T$ or $\R$) each time $l$ intersects a (top or right) side of a square with vertices $(n,m)$, $(n, m+1)$, $(n+1, m)$, and $(n+1,m+1)$ where $n, m \in \mathbb{Z}$, $0 \leq n \leq a$, and $-1 \leq m \leq b - 1$.  
(Each square, along with its interior, covers $(T,p)$ and is triangulated by the lifted arcs of a certain triangulation $\Delta_T$ of $(T,p)$; see left picture of Figure~\ref{R2-to-DZ2}.)
\begin{figure}
\scalebox{0.8}{\includegraphics{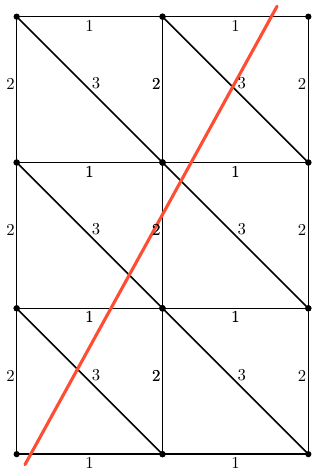}} \quad \raisebox{90pt}{$\longrightarrow$} \quad \scalebox{0.8}{\includegraphics{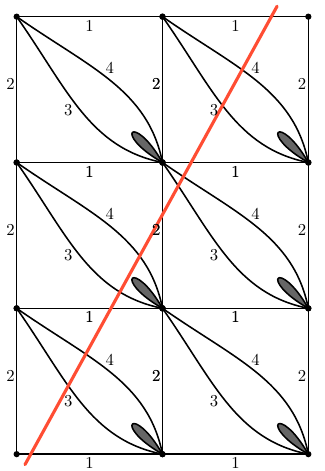}}
\caption{\emph{Left:} A lift of a closed allowable curve in $(T, p)$ along with the lifts of $\Delta_T$ used to calculate the curve's shear coordinates. \emph{Right:} The corresponding closed allowable curve in $(S, p)$ lifted to $\mathbb{R}^2 \setminus D\mathbb{Z}^2$ with the corresponding lifts of the triangulation $\Delta_0$.}
\label{R2-to-DZ2}
\end{figure} 
Viewing each of the integer points in $\mathbb{R}^2$ as a boundary of a ``tiny removed disk,'' $\mathbb{R}^2 \setminus \mathbb{Z}^2$ becomes $\mathbb{R}^2 \setminus D\mathbb{Z}^2$ and each lifted version of $\Delta_T$ becomes a lifted version of $\Delta_0$ (see Figure~\ref{R2-to-DZ2}). 
Each pair of adjacent letters corresponds to the passing of the Type~2 curve's standard lift through three quadrilaterals formed by $\Delta_0$ and so corresponds to a certain contribution to the shear coordinates. 
The sum of these contributions is $[-b, a, b, -a]$. 

Next, suppose $b/a = 0$. This Type~2 curve's arc crossings can be encoded as $\R\R$. 

Now, suppose $-1 \leq b/a < 0$. 
The arc crossings of this Type~2 curve can be encoded in a word formed by recording a letter ($\B$  or $\R$) each time the standard lift intersects the base or right side of a parallelogram with vertices $(n,m)$, $(n+1, m)$, $(n+1, m-1)$, and $(n+2,m-1)$ where $n, m \in \mathbb{Z}$, $-1 \leq n \leq a - 1$, and $1 \leq m \leq b + 1$.  
Again, each pair of adjacent letters corresponds to a certain contribution to the shear coordinates. 

Finally, suppose $-\infty < b/a < -1$.  
If $\mathbb{R}^2 \setminus D\mathbb{Z}^2$ is reflected across the line $y = -x$, each lifted version of $\Delta_0$ is mapped to a (possibly different) lifted version of $\Delta_0$ (with arc indices permuted via $1 \longleftrightarrow 2$ and $3 \longleftrightarrow 4$), a curve crossing a ``Z'' pattern is mapped to a (different) curve crossing a ``backwards Z'' pattern (and vice versa), and the standard lift of a Type~2 curve of slope $d/c$ with $-1 < d/c < 0$ is mapped to a line of slope ${-c}/{-d}$. (This line is a lift of the Type~2 curve of slope ${-c}/{-d}$ since it is isotopic to an integer translate of the Type~2 curve's standard lift.) 
Thus, the shear coordinates of a Type~2 curve of slope $b/a$ with $-\infty < b/a < -1$ can be derived from those of a Type~2 curve of slope $b/a$ with $-1 < b/a < 0$ by interchanging the first and second coordinates, interchanging the third and fourth coordinates, and interchanging the letters $a$ and~$b$. 

\emph{Cases 3 \& 4.} The shear coordinates for a Type~1 curve given by $(b/a, k, \sigma)$ can be calculated using a word determined in the same manner as that of the Type~2 curve of slope $b/a$ (in Case~2) but with one difference, namely, that the Type~1 curve's word may include two new letters ($\C$ and $\s$) that represent the Type~1 curve's standard lift completing a full ($360^\circ$) spiral and completing a partial spiral that makes a nonzero contribution to the shear coordinates, respectively. (Partial spirals that make no contribution to the shear coordinates are not assigned a letter.) 

\emph{Case 5.} This case is done by reflection (analogously to the last subcase of Case~2).

\emph{Cases 6-8.} If $(b/a, k, \sigma) \notin \{ (\infty, 1, -), (0, 1, -), (-1, 0, +) \}$, the shear coordinates can be determined as in Cases 3-5 (with minor adjustments to the argument).
 If $(b/a, k, \sigma) \in \{(\infty, 1, -),  (0, 1, -), (-1, 0, +) \}$, the shear coordinates can be directly calculated from a drawing of the standard lift. 

\end{proof}

\section{Null Tangle Property}

In this section, we prove that the dread torus has the Null Tangle Property and explicitly identify the universal geometric coefficients associated with $\Delta_0$. 

\begin{defn}  [Cf. \protect{\cite[Definitions 7.1 and 7.2]{27}}]
A \emph{tangle} in $(S,p)$ or $(T,p)$ is a finite set of allowable curves, all distinct up to isotopy and not necessarily compatible, with each curve $\lambda$ assigned a \emph{weight} $w_\lambda \in \mathbb{Z}$. 
A tangle is \emph{trivial} if all its weights are $0$. 
The \emph{support} of a tangle is the set of all curves in the tangle with nonzero weights. 
The \emph{shear coordinate vector} of a tangle $\Xi$ with respect to a triangulation $\Delta$ is $\mathbf{b}(\Delta, \Xi) = \sum_{\lambda \in \Xi} w_\lambda \mathbf{b}(\Delta, \lambda)$. 
A tangle $\Xi$ is \emph{null} if $\mathbf{b}(\Delta, \Xi) = \mathbf{0}$ for every triangulation $\Delta$. 
A marked surface is said to have the \emph{Null Tangle Property} if every null tangle is trivial.
\end{defn}

To prove that $(S,p)$ has the Null Tangle Property, we make use of the following three lemmas and three theorems. 
The first lemma was originally proven in the second paragraph of the proof of \cite[Proposition 6.1]{28} and then restated as \cite[Lemma~10.3]{34}. 
The first theorem is a specialization of \cite[Proposition 7.11]{27} to $(S,p)$; the second theorem is \cite[Proposition 4.4]{28}. 
We prove the other lemmas and theorem.

\begin{lem} \label{lem-10.3-in-34}
Let $M$ be a finite subset of $\mathbb{Q} \cup \{ \infty \}$. 
If $f/e \in M$, then there exists $b/a, d/c \in \mathbb{Q} \cup \{ \infty \}$ such that: 
\begin{enumerate}
\item $b/a < d/c < f/e$;
\item $ad-bc=af-be=cf-de=1$; and 
\item there is no $q \in M$ such that $b/a \leq q < f/e$.
\end{enumerate}
\end{lem}

\begin{thm} \label{prop-concl-wt-0}
Let $\Xi$ be a null tangle in $(S,p)$. 
If $(S,p)$ has a triangulation~$\Delta$, an arc $\gamma \in \Delta$, and an allowable curve $\lambda \in \Xi$ such that $b_\gamma(\Delta, \lambda) > 0$ and $b_\gamma(\Delta, \nu) \leq 0$ for every $\nu \in \Xi \setminus \{ \lambda \}$, then $w_\lambda = 0$.
\end{thm}

\begin{lem} \label{support-only-closed}
The support of a null tangle in $(S,p)$ contains only closed curves. 
\end{lem}

\begin{proof}
Suppose $\Xi$ is a null tangle in $(S, p)$. 
By Theorem~\ref{ACs}, each curve in $\Xi$ is either a non-closed curve (i.e., a Type~1 curve) or a closed curve (i.e., a Type~2 or~3 curve). 
Suppose $\Xi$ contains a Type~1 curve of slope $f/e$. 

For each $\nu \in \Xi$ of Type 1~or~2, let $q_\nu$ denote its slope. 
Let $M$ denote the set of all distinct slopes of Type~1 and~2 curves in $\Xi$ (so $f/e \in M$). 
Since $\Xi$ is a tangle, it is finite and so $M$ is finite. 
By Lemma~\ref{lem-10.3-in-34}, there exists $b/a, d/c \in \mathbb{Q} \cup \{ \infty \}$ such that $b/a < d/c < f/e$, $ad-bc=af-be=cf-de=1$, and there is no Type~1 or~2 curve $\nu \in \Xi$ such that $b/a \leq q_\nu < f/e$. 

Since $\Xi$ is finite, it has a curve $\lambda$ that is the greatest in the total order on Type~1 curves of slope $f/e$. 
By Lemma~\ref{Lem-change-spiral-numbers}, we can assume that the spiral number of $\lambda$ is~$1$. (See Remark~\ref{f_n-applied-to-ACs}.) 
Then $\lambda$ is given by either $(f/e, 1, +)$ or $(f/e, 1, -)$. 

\emph{Case 1:  $\lambda$ given by $(f/e, 1, +)$.}
Let $\alpha_1$, $\alpha_2$, $\beta$, and $\gamma$ denote  the arcs given by $(b/a, 0, -)$, $(b/a, 0, +)$, $(d/c, 0, +)$, and $(f/e, 0, +)$, respectively. 
(In Figures~\ref{Support-Case1-1stpic} and~\ref{Support-Case1-2ndpic}, lifts of these arcs are shown in green, red, blue, and purple, respectively.) 
By Theorem~\ref{thm-triangulations}, $\Delta = \{ \alpha_1, \alpha_2, \beta, \gamma \}$ is a triangulation.
From Figure~\ref{Support-Case1-1stpic}, we see that $b_\gamma(\Delta, \lambda) = 1$. 
\begin{figure}
\scalebox{0.175}{\includegraphics{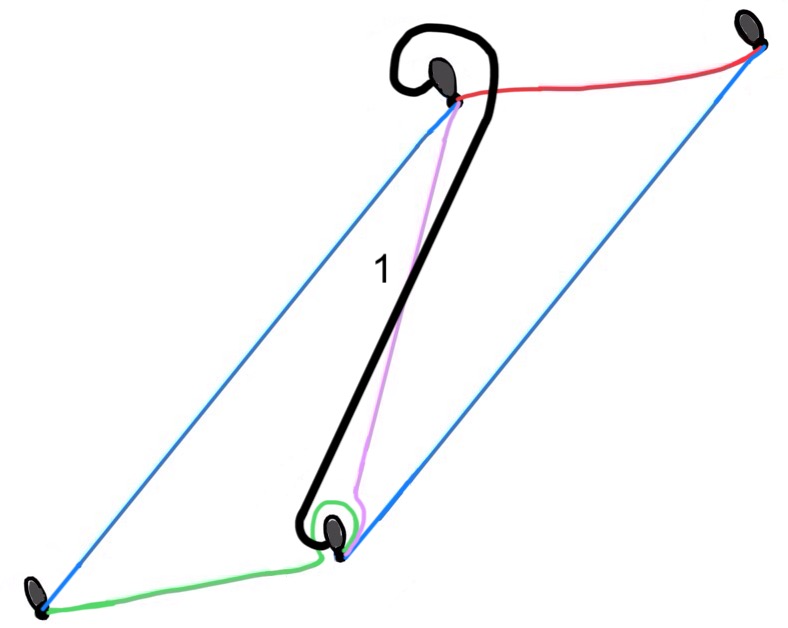}}
\caption{Calculation of $b_\gamma(\Delta, \lambda)$ where $\lambda$ is given by $(f/e, 1, +)$. A lift of $\lambda$ is shown in black.}
\label{Support-Case1-1stpic}
\end{figure}
We also see from this figure that any Type~1 or~2 curve~$\nu$ with $b_\gamma(\Delta, \nu) > 0$ has $b/a < q_\nu \leq f/e$.  
Since there is no Type~1 or~2 curve $\nu \in \Xi$ such that $b/a \leq q_\nu < f/e$, this means that any Type~1 or~2 curve $\nu \in \Xi$ with $b_\gamma(\Delta, \nu) > 0$ has $q_\nu = f/e$.  
From Figure~\ref{Support-Case1-2ndpic}, we see that $b_\gamma(\Delta, \nu) = 0$ if~$\nu$ is the Type~1 curve given by $(f/e, 1, -)$, the Type~2 curve of slope $f/e$, or the Type~3 curve and that $b_\gamma(\Delta, \nu) = -1$ if $\nu$ is a Type~1 curve given by $(f/e, k, \pm)$ where $k < 1$. 
\begin{figure}
\begin{tabular}{cc}
\scalebox{0.175}{\includegraphics{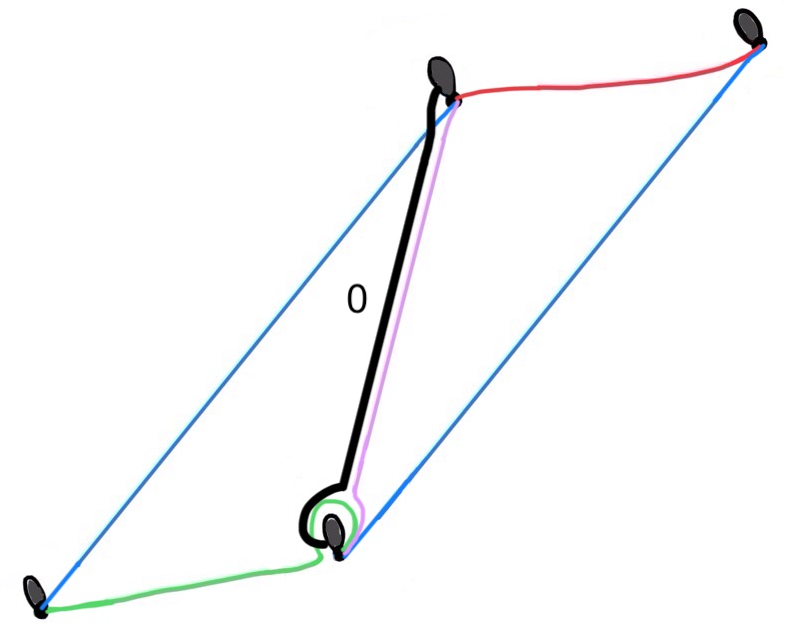}} &
\scalebox{0.175}{\includegraphics{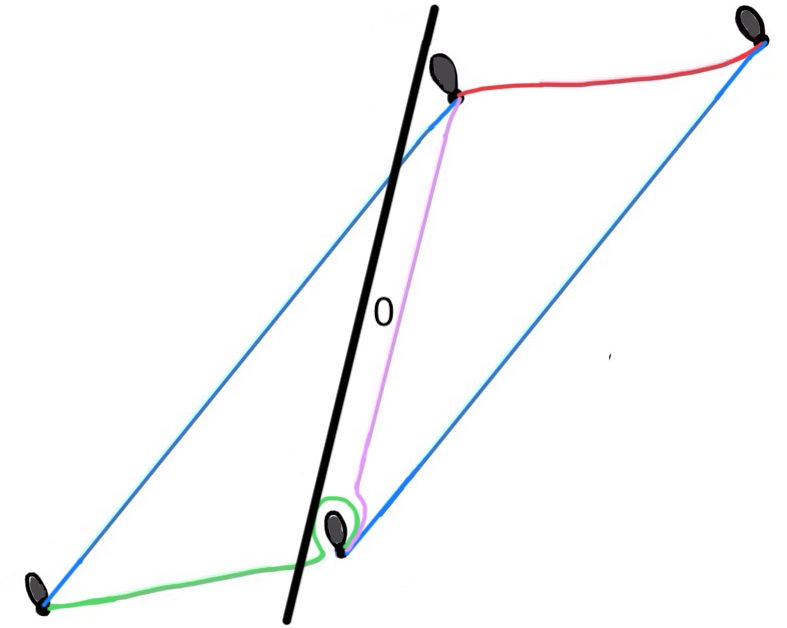}}
\\ 
\scalebox{0.175}{\includegraphics{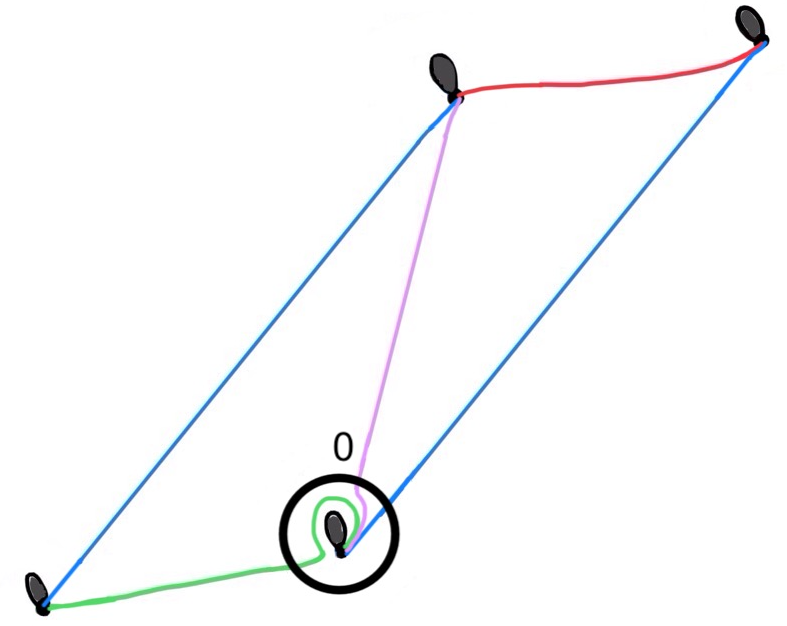}} &
\scalebox{0.175}{\includegraphics{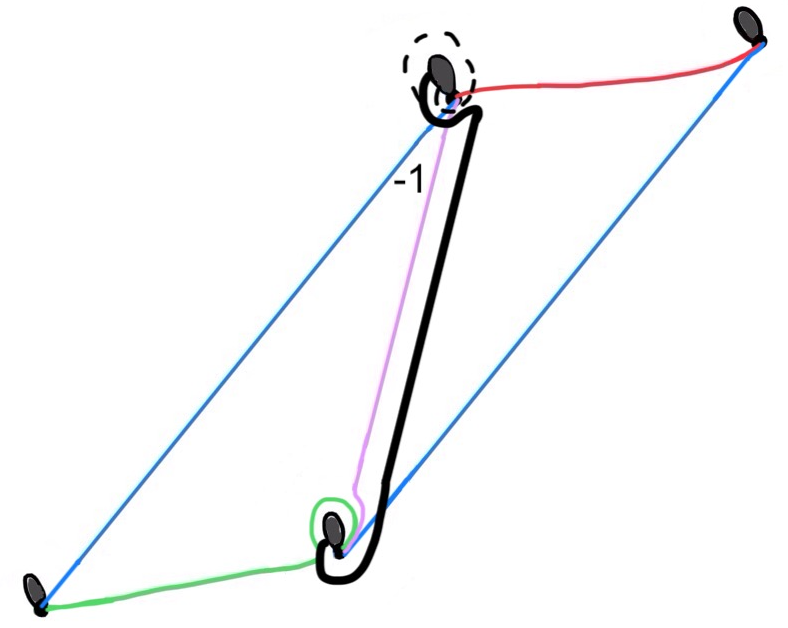}} \\
\end{tabular}
\caption{Calculation of $b_\gamma(\Delta, \nu)$ where $\nu$ is either the Type~1 curve given by $(f/e, 1, -)$ (upper-left image), the Type~2 curve given by $f/e$ (upper-right image), the Type~3 curve (lower-left image), or the Type~1 curve given by $(f/e, 0, \pm)$ (lower-right image). Lifts of $\nu$ are shown in black (with a dashed line for $(f/e, 0, -)$).}
\label{Support-Case1-2ndpic} 
\end{figure} 
Thus, for any $\nu \in \Xi \setminus \{ \lambda \}$, $b_\gamma(\Delta, \nu) \leq 0$.  
By Theorem~\ref{prop-concl-wt-0}, $w_\lambda = 0$. 

\emph{Case 2:  $\lambda$ given by $(f/e, 1, -)$.}
Let $\alpha_1$, $\alpha_2$, $\beta$, and $\gamma$ be the arcs given by $(b/a, -1, +)$, $(b/a, 0, -)$, $(d/c, 0, -)$, and $(f/e, 0, -)$, respectively. 
(In Figures~\ref{Support-Case2-1stpic} and \ref{Support-Case2-2ndpic}, the lifts of these arcs are shown in green, red, blue, and purple, respectively.) 
By Theorem~\ref{thm-triangulations}, $\Delta = \{ \alpha_1, \alpha_2, \beta, \gamma \}$ is a triangulation.
From Figure~\ref{Support-Case2-1stpic}, we see that $b_\gamma(\Delta, \lambda) = 1$ and that any other Type~1 or~2 curve $\nu \in \Xi$ with $b_\gamma(\Delta, \nu) > 0$ has $b/a < q_\nu \leq f/e$. 
\begin{figure}
\scalebox{0.334}{\includegraphics{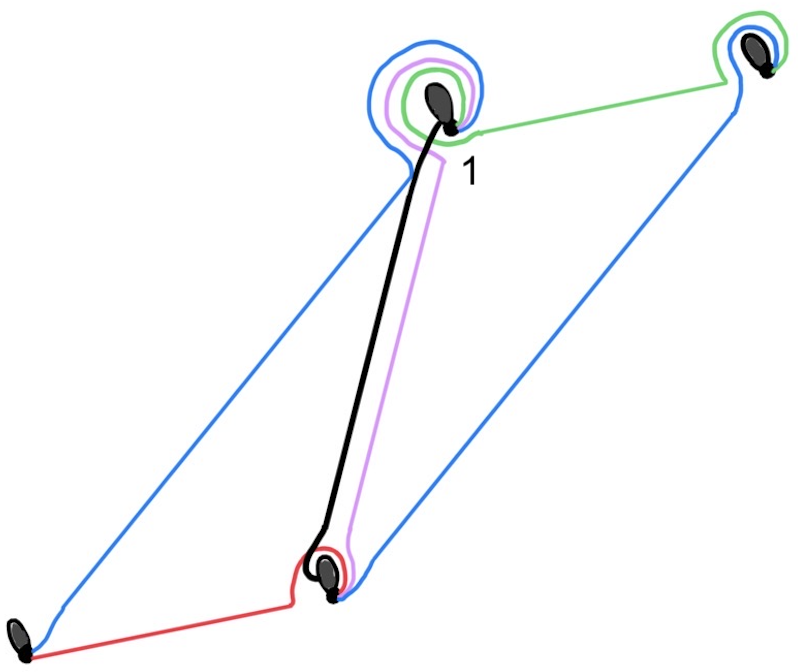}}
\caption{Calculation of $b_\gamma(\Delta, \lambda)$ where $\lambda$ is given by $(f/e, 1, -)$. A lift of $\lambda$ is shown in black.}
\label{Support-Case2-1stpic}
\end{figure}
Since there is no Type~1 or~2 curve $\nu \in \Xi$ such that $b/a \leq q_\nu < f/e$, any Type~1 or~2 curve $\nu \in \Xi$ with $b_\gamma(\Delta, \nu) > 0$ has $q_\nu = f/e$. 
From Figure~\ref{Support-Case2-2ndpic}, we see that $b_\gamma(\Delta, \nu) = 0$ if~$\nu$ is the Type~1 curve given by $(f/e, 0, +)$, the Type~2 curve of slope $f/e$, or the Type~3 curve and that $b_\gamma(\Delta, \nu) = -1$ if $\nu$ is a Type~1 curve given by $(f/e, 0, -)$ or $(f/e, k, \pm)$ where $k < 0$. 
\begin{figure}
\begin{tabular}{cc}
\scalebox{0.141}{\includegraphics{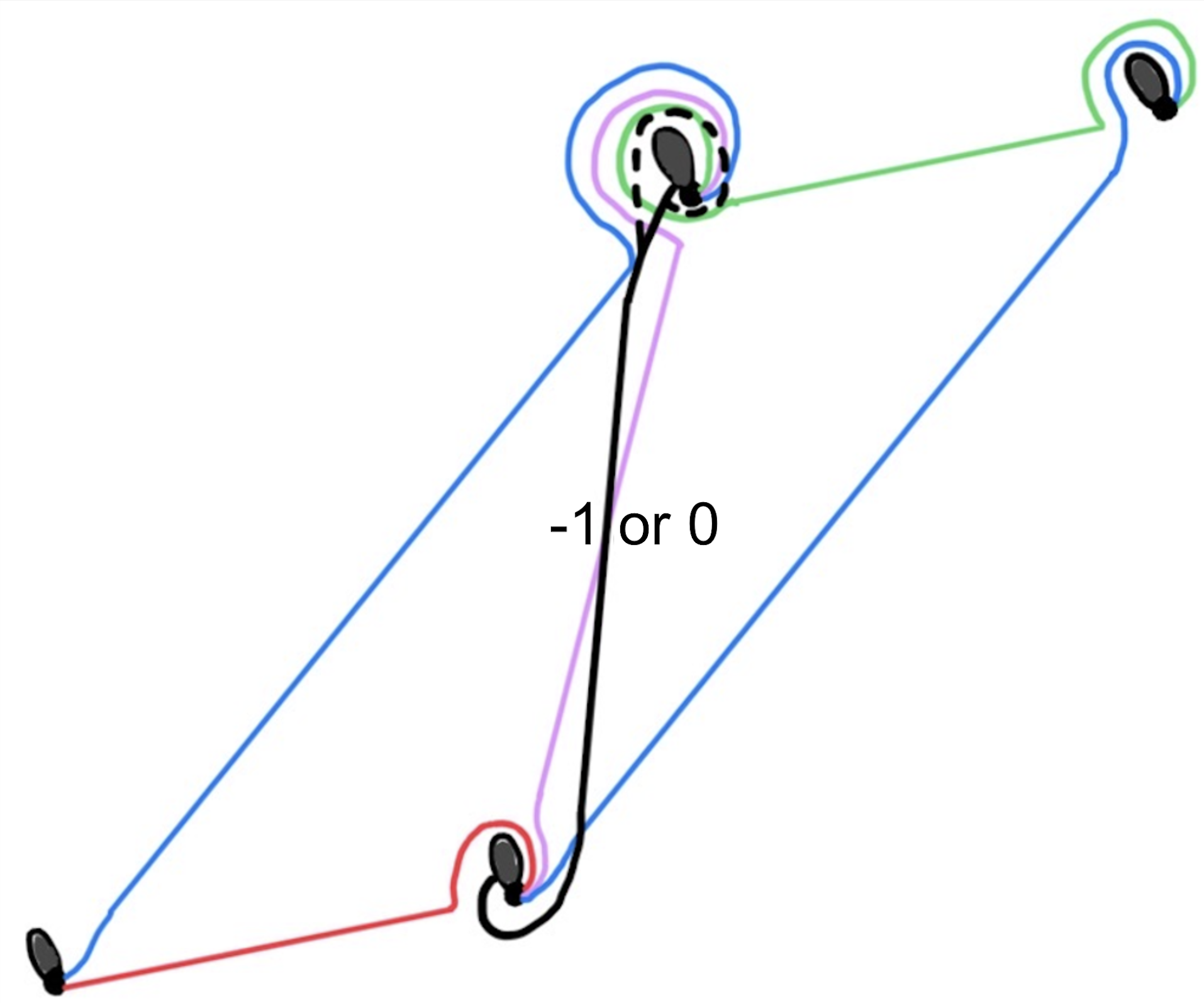}} &
\scalebox{0.187}{\includegraphics{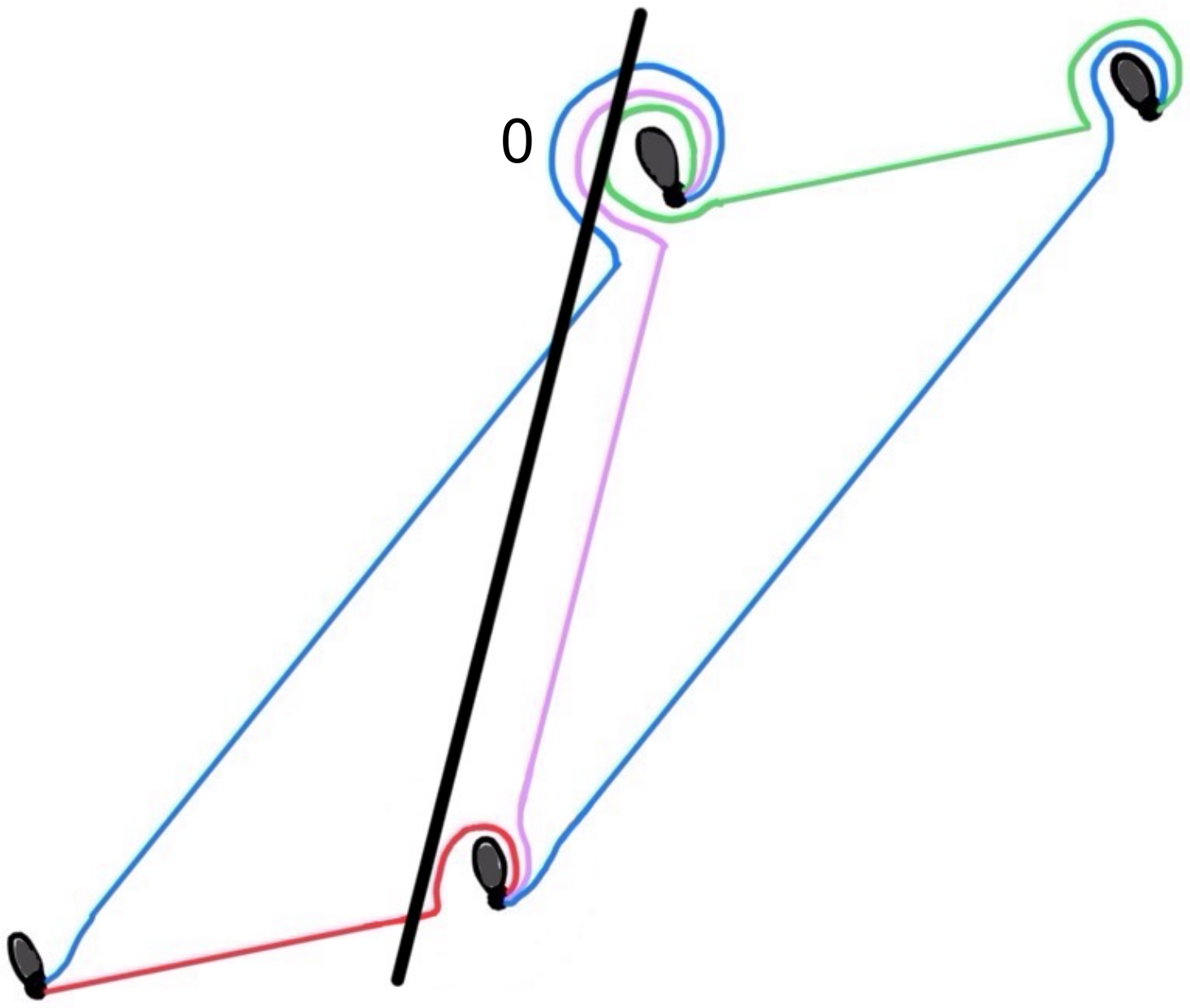}}
\\ 
\scalebox{0.165}{\includegraphics{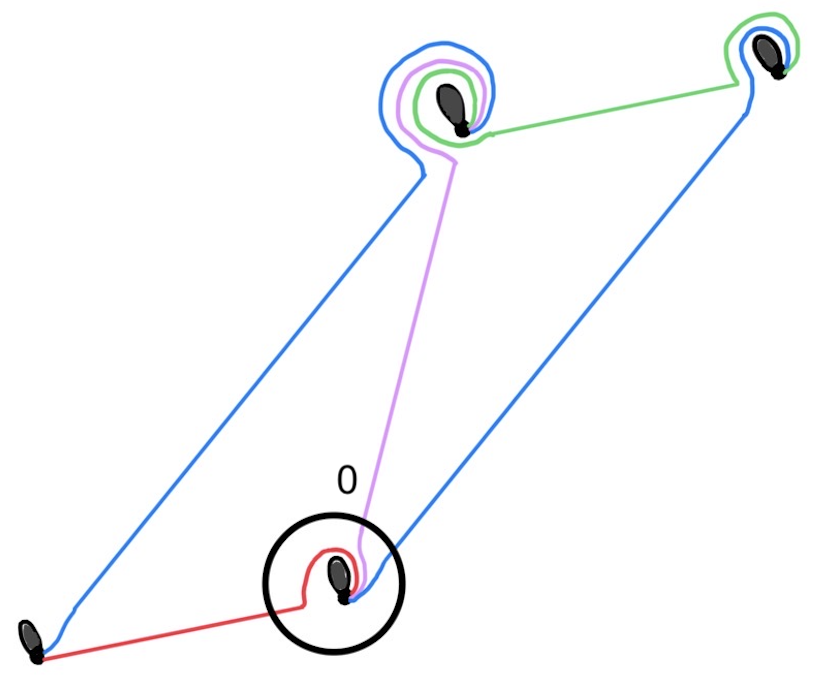}} &
\scalebox{0.1405}{\includegraphics{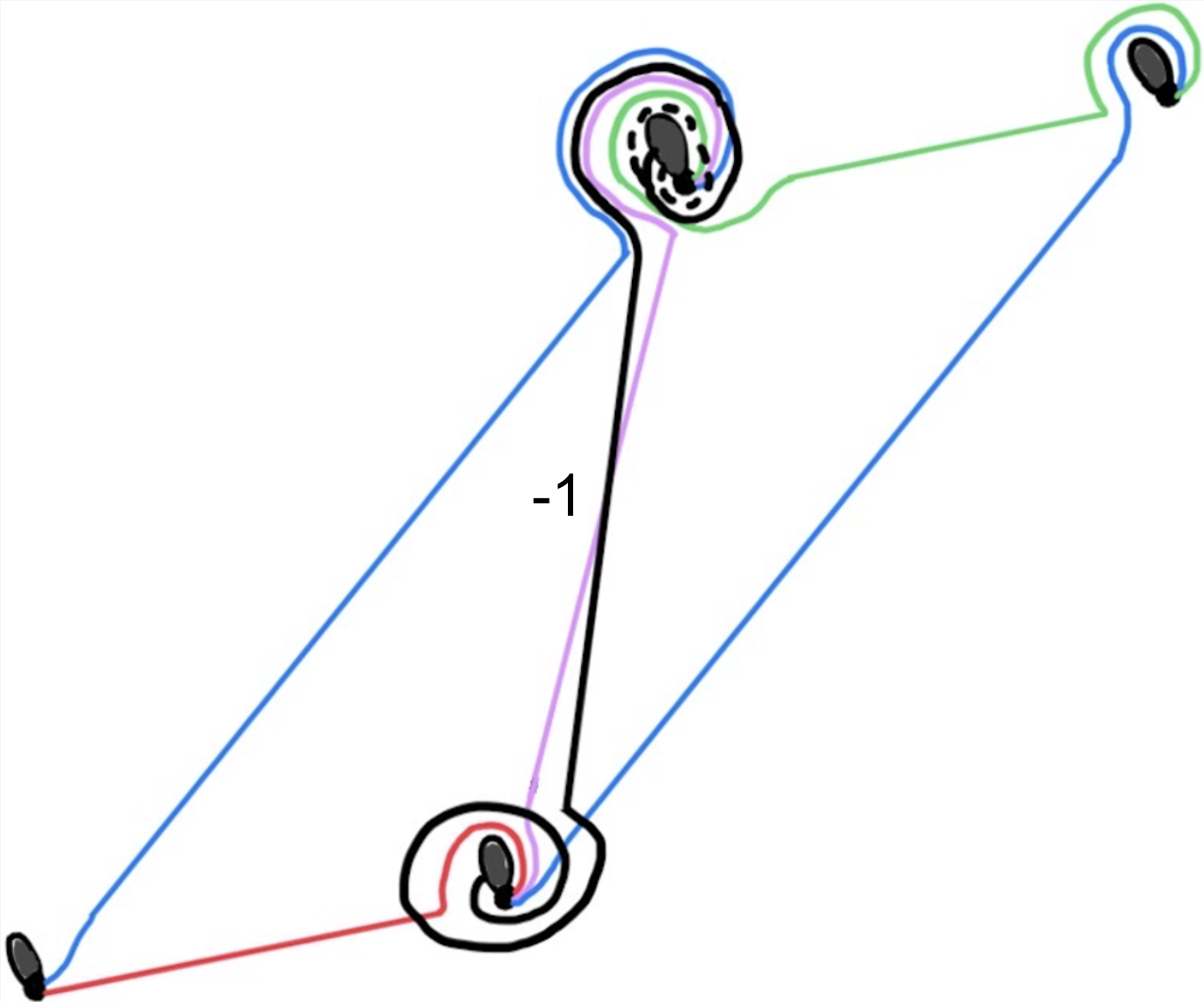}} \\
\end{tabular}
\caption{Calculation of $b_\gamma(\Delta, \nu)$ where $\nu$ is either the Type~1 curve given by $(f/e, 0, \pm)$ (upper-left image), the Type~2 curve given by $f/e$ (upper-right image), the Type~3 curve (lower-left image), or the Type~1 curve given by $(f/e, -1, \pm)$ (lower-right image). Lifts of $\nu$ are shown in black (with dashed lines for $(f/e, 0, -)$ and $(f/e, -1, -)$).}
\label{Support-Case2-2ndpic} 
\end{figure}   
By Theorem~\ref{prop-concl-wt-0}, $w_\lambda = 0$. 

Since, in both cases, $w_\lambda = 0$, $\Xi \setminus \{ \lambda \}$ is a null tangle. 
Thus, if $\Xi \setminus \{ \lambda \}$ contains a Type~1 curve, we can replace $\Xi$ with $\Xi \setminus \{ \lambda \}$ in paragraphs 2-5 of this argument and repeat this part of the argument until there are no Type~1 curves left to remove. 
(This eventually happens since $\Xi$ is finite.) 
This process shows that every Type~1 curve in $\Xi$ has weight $0$ and, therefore, is not in the support of $\Xi$. 
\end{proof}

\begin{thm} \label{art-28-prop-4.4}
The closed allowable curves in $(T,p)$ are the images, under the covering map, of the lines in $\mathbb{R}^2 \setminus \mathbb{Z}^2$ that have rational or infinite slope. 
Any two such allowable curves are isotopic if and only if they lift to lines of the same slope.
\end{thm}

We say that a closed allowable curve in $(T,p)$ that lifts to a line of slope $m$ has \emph{slope} $m$. 
The proofs of the remaining theorem and lemma rely on the close relationship between $(S, p)$ and $(T, p)$, particularly the following three observations: 

\begin{enumerate}
\item \label{component-re-ACs} \emph{There is a bijection from the Type 2 curves of $(S,p)$ to the closed allowable curves of $(T,p)$, induced by the inclusion map $S \hookrightarrow T$, that maps the Type 2 curve of slope $b/a$ to the closed allowable curve of slope $b/a$.} (This follows from Part~\ref{lifts-of-ACs-T2} of Lemma~\ref{lifts-of-ACs}, the fact that each Type~2 curve is uniquely specified by its slope, and Theorem~\ref{art-28-prop-4.4}.) 
\item \label{component-re-arcs} \emph{There is an infinity-to-one, surjective map from the arcs of $(S,p)$ to the arcs of $(T,p)$, induced by the inclusion map $S \hookrightarrow T$,  that maps the arc given by $(b/a, k, \sigma)$ to the arc of slope $b/a$.} (This follows from Theorem~\ref{Prop-4-1-from-28} and Theorem~\ref{Thm-arcs-as-images}.) 
\item \label{component-re-triangulations} \emph{There is an infinity-to-one, surjective map from the triangulations of $(S,p)$ to the triangulations of $(T,p)$ that maps each triangulation $\Delta$ of $(S,p)$ to the (unique) triangulation of $(T,p)$ whose three arcs have the same distinct slopes as the arcs in $\Delta$.} (This follows from Theorem~\ref{thm-triangulations}, \cite[Theorem 4.3]{28}, and the fact that triangulations of $(T,p)$ consist of three arcs.) 
\end{enumerate}

In addition, we frequently use a fact derived from Theorem~\ref{thm-triangulations}, namely, that every triangulation of $(S,p)$ contains exactly two arcs of the same slope. 

\begin{thm} \label{magic-star-doublestar}
Let $\lambda$ be a Type~2 curve in $(S,p)$, and let $\lambda'$ be the closed allowable curve in $(T,p)$ with the same slope as $\lambda$. 
Let $\Delta = \{\alpha, \beta, \gamma_1, \gamma_2 \}$ be a triangulation of $(S,p)$ where $\gamma_1$ and $\gamma_2$ denote the two arcs in $\Delta$ of the same slope, and let $\Delta' = \{\alpha', \beta', \gamma' \}$ be the triangulation of $(T,p)$ where $\alpha'$, $\beta'$, and $\gamma'$ have the same slope, respectively, as $\alpha$, $\beta$, and $\gamma_1$ (equivalently, $\gamma_2$).  Then
$b_\alpha(\lambda, \Delta) = b_{\alpha'}(\lambda', \Delta')$, $b_\beta(\lambda, \Delta) = b_{\beta'}(\lambda', \Delta')$, and $b_{\gamma_1}(\lambda, \Delta) + b_{\gamma_2}(\lambda, \Delta) = b_{\gamma'}(\lambda', \Delta')$.
\end{thm}

\begin{proof}
In what follows, we refer to triangles and quadrilaterals whose sides are arcs in $\Delta'$ or arcs in $\Delta$ and/or the boundary $\delta$ of $S$ by listing their sides in CW order (starting with whichever side is most convenient). 
This naming method uniquely identifies each triangle and quadrilateral except for the triangles whose sides are in~$\Delta'$. 
(Serendipitously, we never have need to uniquely specify one of the triangles whose sides are in $\Delta'$.) 

By Theorem~\ref{thm-triangles-only-slope}, the triangulation $\Delta$ divides $S$ into three triangles: a triangle with sides $\gamma_1$, $\gamma_2$, and $\delta$; a triangle with sides $\alpha$, $\beta$, and $\gamma_1$; and a triangle with sides $\alpha$, $\beta$, and $\gamma_2$. 
By Definition~\ref{Defn-triangulation} and observations (\ref{component-re-arcs}) and (\ref{component-re-triangulations}), this implies that $\Delta'$ divides $T$ into two triangles that each have sides $\alpha'$, $\beta'$, and $\gamma'$. 

Without loss of generality, suppose the triangle with sides $\gamma_1$, $\gamma_2$, and $\delta$ is $\gamma_1 \gamma_2 \delta$ (otherwise, swap $\gamma_1$ and $\gamma_2$ in what follows). 
It is impossible for both $\alpha \beta \gamma_1$ and $\alpha \gamma_2 \beta$ to be triangles or for both $\alpha \gamma_1 \beta$ and $\alpha \beta \gamma_2$ to be triangles, for then $S$ would be a topological disk. 
See Figure~\ref{contradiction-in-magic-double-star}. 
\begin{figure}
\begin{tabular}{cc}
\scalebox{1.6}{\includegraphics{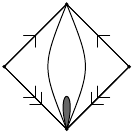}%
\begin{picture}(0,0)(32,-32)
\put(-17,-1){\scriptsize$\gamma_1$}
\put(10,-1){\scriptsize$\gamma_2$}
\put(22,-20){\scriptsize$\alpha$}
\put(20,18){\scriptsize$\beta$}
\put(-27,-20){\scriptsize$\alpha$}
\put(-25,18){\scriptsize$\beta$}
\put(-2,-13){\scriptsize$\delta$}
\end{picture}} & 
\scalebox{1.6}{\includegraphics{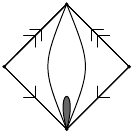}%
\begin{picture}(0,0)(32,-32)
\put(-17,-1){\scriptsize$\gamma_1$}
\put(10,-1){\scriptsize$\gamma_2$}
\put(22,-20){\scriptsize$\beta$}
\put(20,18){\scriptsize$\alpha$}
\put(-27,-20){\scriptsize$\beta$}
\put(-25,18){\scriptsize$\alpha$}
\put(-2,-13){\scriptsize$\delta$}
\end{picture}} \\
\end{tabular}
\caption{\emph{Left:} The surface formed by $\gamma_1 \gamma_2 \delta$, $\alpha \beta \gamma_1$, and $\alpha \gamma_2 \beta$. 
\emph{Right:} The surface formed by $\gamma_1 \gamma_2 \delta$, $\alpha \gamma_1 \beta$, and $\alpha \beta \gamma_2$.} 
\label{contradiction-in-magic-double-star} 
\end{figure} 
(In Figure~\ref{contradiction-in-magic-double-star}, the 
spirals on $\alpha$, $\beta$, $\gamma_1$, and $\gamma_2$ are fixed. 
However, changing these spirals would not alter the  topology.) 
Hence, either $\alpha \beta \gamma_1$ and $\alpha \beta \gamma_2$ are both triangles or $\alpha \gamma_1 \beta$ and $\alpha \gamma_2 \beta$ are both triangles. 
Without loss of generality, suppose $\alpha \beta \gamma_1$ and $\alpha \beta \gamma_2$ are both triangles (otherwise, swap $\alpha$ and $\beta$ in what follows). 
By observations (\ref{component-re-arcs}) and (\ref{component-re-triangulations}), this assumption implies that the two triangles of $\Delta'$ both have the name $\alpha' \beta' \gamma'$. 
The triangulations $\Delta$ and $\Delta'$ are illustrated in Figure~\ref{delta-and-deltaprime-in-magic-double-star} with fixed spirals chosen for the arcs in $\Delta$. 
Since a Type~2 curve does not, up to isotopy, intersect any arc spirals, the choice of spirals does not affect $\mathbf{b}(\Delta, \lambda)$. 
\begin{figure}
\begin{tabular}{cc}
\scalebox{1.6}{\includegraphics{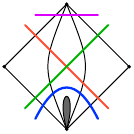}%
\begin{picture}(0,0)(32,-32)
\put(-17,-1){\scriptsize$\gamma_1$}
\put(10,-1){\scriptsize$\gamma_2$}
\put(22,-13){\scriptsize$\beta$}
\put(22,10){\scriptsize$\alpha$}
\put(-27,-13){\scriptsize$\alpha$}
\put(-25,10){\scriptsize$\beta$}
\end{picture}}
\qquad & \qquad
\scalebox{1.6}{\includegraphics{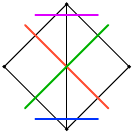}%
\begin{picture}(0,0)(32,-32)
\put(0,-15){\scriptsize$\gamma'$}
\put(22,-13){\scriptsize$\beta'$}
\put(22,10){\scriptsize$\alpha'$}
\put(-29,-15){\scriptsize$\alpha'$}
\put(-26,10){\scriptsize$\beta'$}
\end{picture}}
\end{tabular}
\caption{\emph{Left:} The triangulation $\Delta$ (in black) and the four possibilities for $\lambda$ (in green, orange, blue, and purple). \emph{Right:} The triangulation $\Delta'$ (in black) and the four possibilities for $\lambda'$ (in green, orange, blue, and purple).}
\label{delta-and-deltaprime-in-magic-double-star}
\end{figure} 

Now, we verify that $b_\alpha(\lambda, \Delta) = b_{\alpha'}(\lambda', \Delta')$. 
Since $\alpha'$ forms a side of both triangles named $\alpha' \beta' \gamma'$, it forms the diagonal of a quadrilateral $\beta' \gamma' \beta' \gamma'$. 
Since $\alpha$ forms a side of triangle $\alpha \beta \gamma_1$ and triangle $\alpha \beta \gamma_2$, it forms the diagonal of a quadrilateral $\beta \gamma_1 \beta \gamma_2$. 
Under the inclusion map ${S \hookrightarrow T}$, the quadrilateral $\beta \gamma_1 \beta \gamma_2$ maps to the quadrilateral $\beta' \gamma' \beta' \gamma'$ ($\beta \mapsto \beta'$ and $\gamma_1, \gamma_2 \mapsto \gamma$), $\alpha$ maps to $\alpha'$, and $\lambda$ maps to~$\lambda'$. 
Thus, each time $\lambda'$ crosses $\alpha'$, $\lambda$ crosses $\alpha$ the same way (i.e., either both ``cut a corner'' of their respective quadrilaterals, both cross the diagonal of their respective quadrilaterals in a ``backwards Z'' pattern, or both cross the diagonal of their respective quadrilaterals in a ``Z'' pattern). 
Therefore, each time, the same contribution is made to $b_{\alpha'}(\lambda', \Delta')$ and to $b_\alpha(\lambda, \Delta)$. 

The proof that $b_\beta(\lambda, \Delta) = b_{\beta'}(\lambda', \Delta')$ is analogous to the one used to verify that $b_\alpha(\lambda, \Delta) = b_{\alpha'}(\lambda', \Delta')$ with $\beta'$ forming the diagonal of a quadrilateral $\alpha' \gamma' \alpha' \gamma'$ and~$\beta$ forming the diagonal of a quadrilateral $\alpha \gamma_1 \alpha \gamma_2$. 

Lastly, we verify that $b_{\gamma_1}(\lambda, \Delta) + b_{\gamma_2}(\lambda, \Delta) = b_{\gamma'}(\lambda', \Delta')$. 
Since $\gamma'$ forms a side of both triangles named $\alpha' \beta' \gamma'$, it forms the diagonal of a quadrilateral $\alpha' \beta' \alpha' \beta'$.  
Each time $\lambda'$ crosses $\gamma'$, it enters and exits $\alpha' \beta' \alpha' \beta'$ in one of the four ways shown in the right picture of Figure~\ref{delta-and-deltaprime-in-magic-double-star}. 
Since $\alpha \mapsto \alpha'$, $\beta \mapsto \beta'$, $\gamma_1, \gamma_2 \mapsto \gamma$, and ${\lambda \mapsto \lambda'}$ via the inclusion map $S \hookrightarrow T$, these four possible behaviors of $\lambda'$ correspond to the four possible behaviors of $\lambda$ shown in the left picture of Figure~\ref{delta-and-deltaprime-in-magic-double-star}. 
The corresponding possibilities are shown in the same color.  
From Figure~\ref{delta-and-deltaprime-in-magic-double-star}, we see that the possibilities for $\lambda'$ and $\lambda$ that are shown in green result in contributions of $1$ to $b_{\gamma'}(\lambda', \Delta')$,~$1$ to $b_{\gamma_1}(\lambda, \Delta)$, and $0$ to $b_{\gamma_2}(\lambda, \Delta)$; 
the possibilities for $\lambda'$ and $\lambda$ that are shown in red result in contributions of $-1$ to $b_{\gamma'}(\lambda', \Delta')$, $0$ to $b_{\gamma_1}(\lambda, \Delta)$, and $-1$ to $b_{\gamma_2}(\lambda, \Delta)$; 
the possibilities for $\lambda'$ and $\lambda$ that are shown in blue result in contributions of $0$ to $b_{\gamma'}(\lambda', \Delta')$, $1$ to $b_{\gamma_1}(\lambda, \Delta)$, and $-1$ to $b_{\gamma_2}(\lambda, \Delta)$; and 
the possibilities for $\lambda'$ and $\lambda$ that are shown in purple result in contributions of $0$ to $b_{\gamma'}(\lambda', \Delta')$, $0$ to $b_{\gamma_1}(\lambda, \Delta)$, and $0$ to $b_{\gamma_2}(\lambda, \Delta)$. 
\end{proof}

\begin{lem} \label{lem-Type3-coords}
Let $\{ \alpha, \beta, \gamma_1, \gamma_2 \}$ be a triangulation of $(S,p)$ (where $\gamma_1$ and $\gamma_2$ denote the two arcs of the same slope and $\gamma_1 > \gamma_2$ in the total order on arcs of this shared slope). 
Then the shear coordinate vector of the Type~3 curve with respect to this triangulation is $[0, 0, 1, -1]$. 
\end{lem}

\begin{proof} 
Let $\Delta = \{ \alpha, \beta, \gamma_1, \gamma_2 \}$. 
Figure~\ref{Type3-shear-coords} shows, up to an orientation-preserving invertible linear map on $\mathbb{R}^2$ and a choice of spirals, a lift of the Type~3 curve and the four lifts of $\Delta$ that intersect this lift.  
\begin{figure}
\scalebox{0.2}{\includegraphics{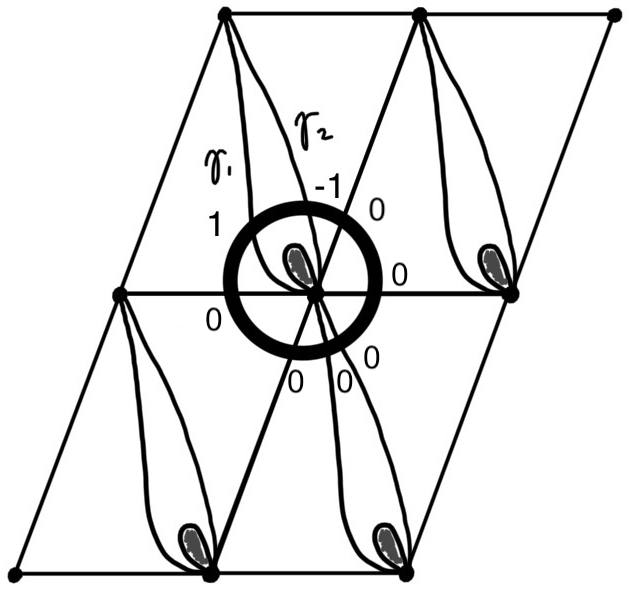}}
\caption{Illustration for the proof of Lemma~\ref{lem-Type3-coords}. 
The thick black circle is the lift of $\lambda$; the unmarked horizontal and diagonal line segments are the lifts of $\alpha$ and $\beta$.}
\label{Type3-shear-coords}
\end{figure}
Since the Type~3 curve does not, up to isotopy, intersect the spirals of any arcs, the choice of spirals does not affect the Type~3 curve's shear coordinates. 
From Figure~\ref{Type3-shear-coords}, we calculate that  the Type~3 curve's shear coordinate vector with respect to $\Delta$ is $[0, 0, 1, -1]$. 
\end{proof}

\begin{thm} \label{Thm-NTP}
The dread torus has the Null Tangle Property.
\end{thm}

\begin{proof}
Suppose $\Xi$ is a null tangle in $(S,p)$. 
By Lemma~\ref{support-only-closed}, the curves in $\Xi$ with nonzero weight are all of Type~2 or~3. 
Let $\lambda$ be the Type~3 curve, and let $\bar{\Xi}$ be the tangle $\Xi \setminus \{ \lambda \}$. 
Then, for any triangulation $\Delta$ of $(S,p)$, 
\begin{equation} \label{1st-vector-eqn}
u \mathbf{b} (\lambda, \Delta) + \sum_{ \nu \in \bar{\Xi}} w_\nu \mathbf{b} (\nu, \Delta) = \mathbf{0}
\end{equation}
where $u = w_\lambda$ if $\lambda \in \Xi$ and $u=0$ otherwise. 
Let $\alpha$, $\beta$, $\gamma_1$, and $\gamma_2$ be the arcs of $\Delta$ where $\gamma_1$ and $\gamma_2$ denote the two arcs of the same slope and $\gamma_1 > \gamma_2$ in the total order on arcs of this shared slope. 
By Lemma~\ref{lem-Type3-coords}, $\mathbf{b} (\lambda, \Delta) = [0, 0, 1, -1]$.
Substituting this vector into equation~(\ref{1st-vector-eqn}) and then breaking the resulting vector equation into components yields: 
$$\sum_{ \nu \in \bar{\Xi}} w_\nu b_\alpha (\nu, \Delta) = 0, \quad 
\sum_{ \nu \in \bar{\Xi}} w_\nu b_\beta (\nu, \Delta) = 0, $$
$$u + \sum_{ \nu \in \bar{\Xi}} w_\nu b_{\gamma_1} (\nu, \Delta) = 0, \quad \text{and} \quad {-u} + \sum_{ \nu \in \bar{\Xi}} w_\nu b_{\gamma_2} (\nu, \Delta) = 0.$$  
Thus: 
$$\sum_{ \nu \in \bar{\Xi}} w_\nu \big( b_{\gamma_1} (\nu, \Delta) + b_{\gamma_2} (\nu, \Delta) \big ) = \sum_{ \nu \in \bar{\Xi}} w_\nu b_{\gamma_1} (\nu, \Delta) + \sum_{ \nu \in \bar{\Xi}} w_\nu b_{\gamma_2} (\nu, \Delta) = 0.$$ 

Let $\alpha'$, $\beta'$, and $\gamma'$ be the arcs of $(T,p)$ that have the same slopes as $\alpha$, $\beta$, and~$\gamma_1$, respectively. 
Let $\Delta' = \{ \alpha', \beta', \gamma' \}$. 
By observation (\ref{component-re-triangulations}), $\Delta'$ is a triangulation of $(T,p)$. 
By observation (\ref{component-re-ACs}), there is a tangle $\bar{\Xi}'$ in $(T,p)$ that consists of the closed allowable curves with the same slopes and same weights as the Type~2 curves in $\bar{\Xi}$. 
By Theorem~\ref{magic-star-doublestar} and observation (\ref{component-re-ACs}), for any Type~2 $\nu \in \bar{\Xi}$ and $\nu' \in \bar{\Xi'}$ of the same slope, $b_\alpha (\nu, \Delta) = b_{\alpha'} (\nu', \Delta')$, $b_\beta (\nu, \Delta) = b_{\beta'} (\nu', \Delta')$, and $b_{\gamma_1}(\nu, \Delta) + b_{\gamma_2}(\nu, \Delta) = b_{\gamma'}(\nu', \Delta')$. 
Thus: 
$$\sum_{ \nu' \in \bar{\Xi}'} w_\nu b_{\alpha'} (\nu', \Delta') = 0, \quad \sum_{ \nu' \in \bar{\Xi}'} w_\nu b_{\beta'} (\nu', \Delta') = 0, \quad \text{and} \quad \sum_{ \nu' \in \bar{\Xi}'} w_\nu b_{\gamma'} (\nu', \Delta') = 0.$$ 
Together, these equations are equivalent to vector equation $\sum_{ \nu' \in \bar{\Xi}'} w_\nu \mathbf{b} (\nu', \Delta') = \mathbf{0}$. 
Letting $\Delta$ vary over all possible triangulations of $(S,p)$, the triangulation $\Delta'$ also varies over all possible triangulations of $(T,p)$.  
Thus, this vector equation implies that $\bar{\Xi}'$ is a null tangle. 
Since $(T,p)$ possesses the Null Tangle Property (Section~\ref{Sect-Intro}), this means that the weights of all the curves in $\bar{\Xi}'$ are~$0$. 
Substituting these weights into equation (\ref{1st-vector-eqn}) yields $u \mathbf{b} (\lambda, \Delta)  = \mathbf{0}$.  Since $\mathbf{b} (\lambda, \Delta) \neq \mathbf{0}$, $u = 0$. 
Therefore, $w_\nu = 0$ for all $\nu \in \Xi$.
\end{proof}

According to \cite[Theorems 2.1 and 7.3]{27}, if an unpunctured marked surface has the Null Tangle Property, universal geometric coefficients exist for every initial exchange matrix specified by a triangulation of that surface. 
We now explicitly identify the universal geometric coefficients for $\Delta_0$ and explain how universal geometric coefficients can be obtained for the other  triangulations of $(S,p)$. 

\begin{defn}[\protect{\cite[Definitions 2.1 and 2.2]{26}}]  \label{Defn-tropical-semifield}
Let $R$ be either $\mathbb{Z}$ or a subfield of $\mathbb{R}$ that contains $\mathbb{Q}$ as a subfield. 
Let $I$ be a (possibly uncountable) indexing set, $u_i$ a formal symbol for each $i \in I$ (where $u_i \neq u_j$ if $i \neq j$), and $a_i \in R$ for each $i \in I$. 
A \emph{tropical semifield over $R$} is a set of all formal products of the form $\prod_{i \in I} u_i^{a_i}$ endowed with a multiplication $\cdot$ and an addition $\oplus$ defined by
$$\prod_{i \in I} u_i^{a_i} \cdot \prod_{i \in I} u_i^{b_i} = \prod_{i \in I} u_i^{a_i + b_i} \quad \text{and} \quad \prod_{i \in I} u_i^{a_i} \oplus \prod_{i \in I} u_i^{b_i} = \prod_{i \in I} u_i^{\min(a_i, b_i)}.$$ 
\end{defn}

Consider a geometric cluster algebra that has a set of ``combinatorial data'' consisting of an exchange matrix $B$ and coefficients $y_1, y_2, ..., y_n$. 
Using the notation of Definition~\ref{Defn-tropical-semifield}, this implies that, for each $j \in \{ 1, 2, ..., n \}$, $y_j = \prod_{i \in I} u_i^{a_{ij}}$ for some $a_{ij} \in R$. 
Hence, for each $j \in \{ 1, 2, ..., n \}$, $y_j$ can be uniquely identified by the ordered list $(a_{ij} : i \in I)$ (and, conversely, this ordered list uniquely identifies $y_j$). 
This ordered list is typically thought of as a column (of possibly infinite height) and appended below the $j$th column of $B$ \cite{26, 27}. 
The matrix $B$ with these appended columns is an $(n + |I|) \times n$ matrix denoted $\tilde{B}$ \cite{26}.

\begin{defn}[\protect{\cite[Definition 2.3]{26}}] \label{Defn-extended-exchange-matrix}
The matrix $\tilde{B}$ described above is called the \emph{extended exchange matrix} for $B$. 
Its last $|I|$ rows are called its \emph{coefficient rows}. 
\end{defn}

If $B$ is specified by a triangulation $\Delta$, it can be calculated as a \emph{signed adjacency matrix} \cite{CATS1}. 

\begin{defn}[Cf. \protect{\cite[Definition 4.1]{CATS1}}] 
Let $\Delta = \{ \alpha_1, \alpha_2, \alpha_3, \alpha_4 \}$ be a triangulation of $(S, p)$. 
The \emph{signed adjacency matrix} $B(\Delta)$ is a $4 \times 4$ matrix indexed by the arcs of~$\Delta$.
Each entry $b_{\alpha_i \alpha_j}$ is the sum, over the three triangles of $\Delta$, of a 1 if $\alpha_i$ immediately precedes $\alpha_j$ in a CW path around a given triangle, a -1 if $\alpha_j$ immediately precedes $\alpha_i$ in a CW path around a given triangle, and a $0$ otherwise.
\end{defn}

For example, we see from the right picture of Figure~\ref{FigforTriangulations} that
$$B(\Delta_0) =
\begin{bsmallmatrix*}[r]
    0 & 2 & -1 & -1 \\
    -2 & 0 & 1 & 1 \\
    1 & -1 & 0 & 1 \\
    1 & -1 & -1 & 0 \\
\end{bsmallmatrix*}.$$
Let $\tilde{B}(\Delta_0)$ denote the extended exchange matrix formed by appending the shear coordinate vectors listed in Theorem~\ref{Thm-shear-coords-of-ACs} to $B(\Delta_0)$ as coefficient rows. 
The following corollary is a consequence of Theorem~\ref{Thm-NTP} and \cite[Theorems 2.1 and 7.3]{27}. 

\begin{cor} \label{Cor-universal-geom-coeffs-wrt-Delta0}
The coefficient rows of $\tilde{B}(\Delta_0)$ 
represent geometric coefficients that are universal (over $\mathbb{Z}$ and $\mathbb{Q}$) for $B(\Delta_0)$. 
\end{cor}

\begin{rem} \label{Rem-universal-geom-coeffs-wrt-notDelta0}
Matrix mutations can be naturally extended to apply to extended exchange matrices (cf. \cite{CA4, 26}). 
By \cite[Theorems 2.1, 3.9, and 7.3]{27}, the result of applying a matrix mutation $\mu_k$ to $\tilde{B}(\Delta_0)$ is the extended exchange matrix whose upper $n \times n$ submatrix equals $B(\Delta^k)$
(where $\Delta^k$ is the triangulation of $(S,p)$ derived from $\Delta_0$ by a process called \emph{flipping the $k$th arc of $\Delta_0$}; see \cite{CATS1, CATS2} for more details) 
and whose coefficient rows represent geometric coefficients that are universal (over $\mathbb{Z}$ and $\mathbb{Q}$) for $B(\Delta^k)$. 
Any triangulation $\Delta$ of $(S,p)$ can be obtained from $\Delta_0$ by a finite sequence of arc flips (cf. \cite{CATS1}), so the universal geometric coefficients for $B(\Delta)$ can be obtained from $\tilde{B}(\Delta_0)$ by a finite sequence of matrix mutations. 
\end{rem}

\section{Mutation Fan}
 
In this section, we describe the mutation fan of the dread torus.
We assume the reader is familiar with common terminology and facts regarding fans. 
We start by stating some definitions and facts specifically pertaining to mutation fans.

The \emph{mutation fan} $\mathcal{F}_{B}$ of an $n \times n$ exchange matrix $B$ is a fan in $\mathbb{R}^n$ that serves as a geometric description of how repeated mutations of an extended exchange matrix for $B$ affect its coefficient rows. 
(The precise definition of $\mathcal{F}_{B}$ is not needed for this investigation but can be found in  \cite[Definition 5.12] {26}.) 
The mutation fan $\mathcal{F}_B$ is complete \cite[Theorem 5.13]{26}. If $\Delta$ and $\Delta'$ are triangulations of the same marked surface, then $\mathcal{F}_{B(\Delta)}$ and $\mathcal{F}_{B(\Delta')}$ are isomorphic (via the map described in \cite[Definition 4.1]{27}; see \cite[Proposition 7.3]{27}). 

\begin{defn}[Cf. \protect{\cite[Definition 4.8]{27}}] \label{Defn-rational-quasi-lamination-fan} 
Let $\Lambda$ be a collection of pairwise compatible allowable curves in $(S,p)$. 
Let $\Delta$ be a triangulation of $(S,p)$, and 
let $C_\Lambda$ be the nonnegative span of $\{ \mathbf{b}(\Delta, \lambda) : \lambda \in \Lambda \}$. 
The \emph{rational quasi-lamination fan} for~$\Delta$, denoted $\mathcal{F}_{\mathbb{Q}}(\Delta)$, is the collection of all possible $C_\Lambda$. 
\end{defn}

As suggested by the name, $\mathcal{F}_{\mathbb{Q}}(\Delta)$ is a rational fan \cite[Theorem 4.10]{27}. 
Furthermore, for each collection $\Lambda$, \cite[Theorem 4.4]{27} implies that the shear coordinate vectors of the curves in $\Lambda$ are linearly independent. 
Thus, the dimension of each $C_\Lambda$ is $|\Lambda|$ and the faces of each $C_\Lambda$ are the nonnegative spans of the shear coordinate vectors of the curves in every subset of $\Lambda$. 
The following theorem describes how $\mathcal{F}_\mathbb{Q}(\Delta)$ relates to $\mathcal{F}_{B(\Delta)}$. 
(The first two items are a combination of \cite[Theorem 4.10]{27}, \cite[Corollary~7.14]{27}, and \cite[Proposition 6.10]{26}; the third item is part of \cite[Theorem 3.10]{34}.) 

\begin{thm} \label{Thm-conglomerate1}
For a marked surface with the Null Tangle Property:
\begin{itemize} 
\item All rational cones of $\mathcal{F}_{B(\Delta)}$ are in $\mathcal{F}_{\mathbb{Q}}(\Delta)$. 
\item The full-dimensional cones of $\mathcal{F}_{B(\Delta)}$ and $\mathcal{F}_{\mathbb{Q}}(\Delta)$ are the same. 
\item All codimension-1 cones of $\mathcal{F}_{\mathbb{Q}}(\Delta)$ are in $\mathcal{F}_{B(\Delta)}$. 
\end{itemize}
\end{thm}

We now turn to describing $\mathcal{F}_{B(\Delta_0)}$. 
Since $|\Delta_0| = 4$ and $\mathcal{F}_{B(\Delta_0)}$ is complete, the cones of $\mathcal{F}_{B(\Delta_0)}$ fill $\mathbb{R}^4$. 
By Theorems~\ref{Thm-NTP} and~\ref{Thm-conglomerate1}, all the rational cones of $\mathcal{F}_{B(\Delta_0)}$ are in $\mathcal{F}_\mathbb{Q}(\Delta_0)$. 
Theorems~\ref{Thm-shear-coords-of-ACs} and~\ref{Thm-maxlcolsofACs} allow us, in principle, to construct $\mathcal{F}_\mathbb{Q}(\Delta_0)$. 
The maximal cones of $\mathcal{F}_\mathbb{Q}(\Delta_0)$ are the cones formed from the maximal collections of pairwise compatible allowable curves. 
By Theorem~\ref{Thm-maxlcolsofACs}, these collections come in three types (Types~A,~B, and~C). 
(We reuse the names Type~A, Type~B, and Type~C for the cones formed from these collections and reuse the names Type~1, Type~2, and Type~3 from Theorem~\ref{ACs} for the 1-dimensional cones, i.e., \emph{rays}, formed from singleton sets of allowable curves.  
We also say that a Type~1 ray is given by $(m, k, \sigma)$ if its Type~1 curve is given by $(m, k, \sigma)$.) 
Since the dimension of $C_\Lambda$ is $|\Lambda|$,  the Type~A cones are 4-dimensional, the Type~B cones are 3-dimensional, and the Type~C cones are 2-dimensional (see Remark~\ref{Rem-cardinality-of-maxl-cols}). 
Hence, by Theorems~\ref{Thm-NTP} and~\ref{Thm-conglomerate1}, the Type~A and~B cones, as well as their faces, are in $\mathcal{F}_{B(\Delta_0)}$ but the Type~C cones and the Type~3 ray 
might not be in $\mathcal{F}_{B(\Delta_0)}$. 
(Later, we conjecture that the Type~C cones and the Type~3 ray actually belong to $\mathcal{F}_{B(\Delta_0)}$; see Conjecture~\ref{conj-all-in-mut}.) 

Since the Type~A cones are 4-dimensional, it is impossible to draw $\mathcal{F}_\mathbb{Q}(\Delta_0)$ 
on 2-dimensional paper. 
However, we will give some illustrations of the geometry and combinatorics of $\mathcal{F}_\mathbb{Q}(\Delta_0)$. 

The sum of the coordinates of each shear coordinate vector listed in Theorem~\ref{Thm-shear-coords-of-ACs} is equal to either $0$, $1$, or $-1$.  
Let $\mathcal{S}_0$ denote the linear 3-space consisting of all vectors whose coordinate sum is~$0$, and let $\mathcal{S}_1$ and $\mathcal{S}_{-1}$ denote the affine 3-spaces consisting of all vectors whose coordinate sum is~$1$ and all vectors whose coordinate sum is~$-1$, respectively. 
The rays in $\mathcal{F}_\mathbb{Q}(\Delta_0)$ that correspond to shear coordinate vectors whose sum is $0$ lie entirely in $\mathcal{S}_0$, the rays in $\mathcal{F}_\mathbb{Q}(\Delta_0)$ that correspond to shear coordinate vectors whose sum is $1$ intersect (in exactly one point) $\mathcal{S}_1$, and the rays in $\mathcal{F}_\mathbb{Q}(\Delta_0)$ that correspond to shear coordinate vectors whose sum is $-1$ intersect (in exactly one point) $\mathcal{S}_{-1}$. 
In particular, the Type~3 ray, the Type~2 rays, and the Type~1 rays given by either $(m, 1, -)$ where $m > -1$ or $(m, 0, +)$ where $m < -1$ all lie in $\mathcal{S}_0$; 
the Type~1 rays given by either $(m, k, -)$ where $m > -1$ and $k > 1$, $(m, k, +)$ where $m > -1$ and $k \geq 1$, or $(m, k, \pm)$ where $m \leq -1$ and $k > 0$ all intersect $\mathcal{S}_1$; and 
the Type~1 rays given by either $(m, k, \pm)$ where $m \geq -1$ and $k < 1$, $(m, k, -)$ where $m < -1$ and $k \leq 0$, or $(m, k, +)$ where $m < -1$ and $k < 0$ all intersect $\mathcal{S}_{-1}$. 
Hence, by Theorems~\ref{Thm-shear-coords-of-ACs} and~\ref{Thm-maxlcolsofACs}, each cone of $\mathcal{F}_\mathbb{Q}(\Delta_0)$ is contained in one of the two halfspaces bounded by $\mathcal{S}_0$ except for those having the nonnegative span of the shear coordinate vectors of the Type~1 curves given by $(-1, 1, -)$ and $(-1, 0, +)$ as a face.

\begin{figure}
\scalebox{0.42}{\includegraphics{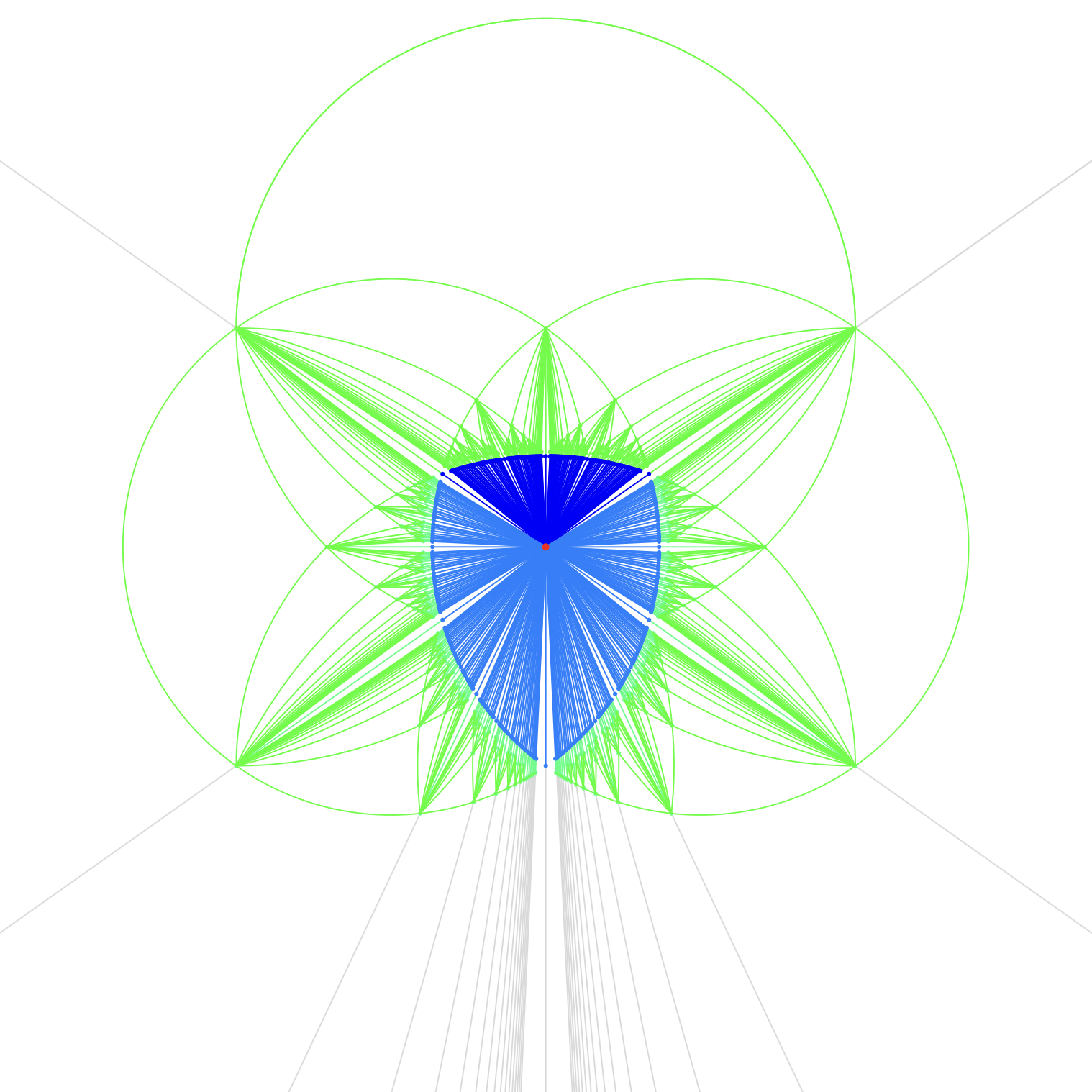}}
\caption{Stereographic projection of the intersection of the sphere $O$ with the cones of $\mathcal{F}_{\mathbb{Q}} \big{|} \mathcal{S}_0$ (through the point on $O$ antipodal to the intersection point of $O$ and the Type~3 ray).}
\label{SUM=0}
\end{figure}

\begin{figure}
\scalebox{0.42}{\includegraphics{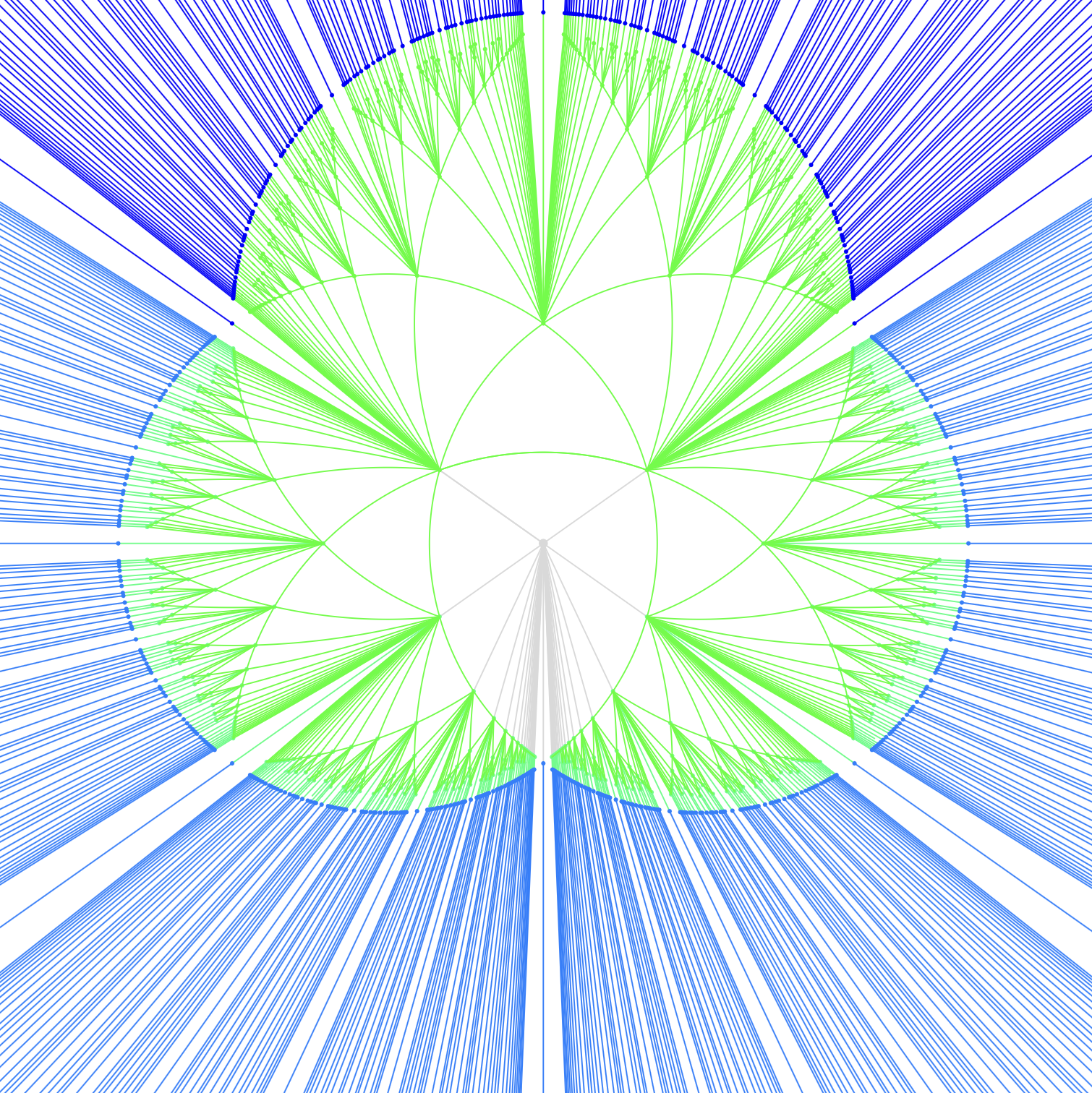}}
\caption{Stereographic projection of intersection of the sphere $O$ with the cones of $\mathcal{F}_{\mathbb{Q}} \big{|} \mathcal{S}_0$ (through the intersection point of $O$ and the Type~3 ray).}
\label{SUM=0reversed}
\end{figure}

Since $\mathcal{S}_0$ is a linear 3-space and the cones of $\mathcal{F}_\mathbb{Q}(\Delta_0)$ are the nonnegative spans of finite sets of vectors, the set $\mathcal{F}_{\mathbb{Q}} \big{|} \mathcal{S}_0 = \{ C \cap \mathcal{S}_0 : C \in \mathcal{F}_\mathbb{Q}(\Delta_0) \}$ forms a fan (cf. \cite[p.~28 and Section~7.1]{Ziegler}) whose cones all have dimension at most~$3$. 
Thus, $\mathcal{F}_{\mathbb{Q}} \big{|} \mathcal{S}_0$ can be graphed on (2-dimensional) paper in a way that is both geometrically accurate and combinatorially meaningful, namely, by first intersecting its cones with the unit sphere $O$ in $\mathcal{S}_0$ that is centered at the origin and then stereographically projecting these intersections through a suitable point on $O$. 
Figures \ref{SUM=0} and \ref{SUM=0reversed} show two different choices of stereographic projection. 
The cones of $\mathcal{F}_{\mathbb{Q}} \big{|} \mathcal{S}_0$ that belong to $\mathcal{F}_\mathbb{Q}(\Delta_0)$ form a subfan of $\mathcal{F}_{\mathbb{Q}} \big{|} \mathcal{S}_0$. 

In Figures~\ref{SUM=0} and~\ref{SUM=0reversed}, the rays of $\mathcal{F}_{\mathbb{Q}} \big{|} \mathcal{S}_0$ have become dots, the 2-dimensional cones of $\mathcal{F}_{\mathbb{Q}} \big{|} \mathcal{S}_0$ have become circular arcs, and the 3-dimensional cones of $\mathcal{F}_{\mathbb{Q}} \big{|} \mathcal{S}_0$ have become curvilinear triangles. 
Colorful dots, arcs, and triangles in both figures correspond to cones of $\mathcal{F}_{\mathbb{Q}} \big{|} \mathcal{S}_0$ that are also in $\mathcal{F}_\mathbb{Q}(\Delta_0)$; gray ones correspond to cones of $\mathcal{F}_{\mathbb{Q}} \big{|} \mathcal{S}_0$ that are intersections of cones of $\mathcal{F}_\mathbb{Q}(\Delta_0)$ with $\mathcal{S}_0$ but are not in $\mathcal{F}_\mathbb{Q}(\Delta_0)$. 
In particular, the red dot (in the center of Figure \ref{SUM=0}; off the page at the ``point at infinity'' in Figure~\ref{SUM=0reversed}) represents the Type~3 ray. 
The blue dots (forming the outline of the ``shield'' shape in Figure~\ref{SUM=0} and the outline of the ``Venn diagram'' shape in Figure~\ref{SUM=0reversed}) represent the Type~2 rays (dark blue for Type~2 curves with nonnegative slope; medium blue for Type~2 curves with negative slope). 
(In both figures, the dark blue dot farthest to the left corresponds to slope $\infty$ and the slopes corresponding to the blue dots decrease as one travels CW from that dark blue dot around the [red or gray] dot in the center of the picture.) 
The green dots represent the Type~1 rays given by either $(m, 1, -)$ with $m > -1$ or $(m, 0, +)$ with $m < -1$. 
The gray dot (off the page at the ``point at infinity'' in Figure~\ref{SUM=0}; in the center of Figure \ref{SUM=0reversed}) represents the 1-dimensional intersection of the 2-dimensional cone of $\mathcal{F}_\mathbb{Q}(\Delta_0)$ consisting of the nonnegative span of the shear coordinate vectors of the Type~1 curves given by $(-1, 1, -)$ and $(-1, 0, +)$. 

Observe that, in both figures, the red dot is joined to each blue dot by a blue arc. 
These blue arcs represent the Type~C cones.  
In addition, each blue dot\textemdash \emph{except} for the medium blue dot that is centered horizontally in both pictures\textemdash is joined to a green dot by a green arc. 
In particular, although it cannot be seen from either figure, a blue dot and a green dot are joined by an arc if and only if they represent a Type~2 ray and a Type~1 ray whose curves have the same slope (cf. Theorem~\ref{Thm-maxlcolsofACs}). 
Each of these green arcs represents a 2-dimensional cone that is a face of exactly two Type~B cones (cf. Theorem~\ref{Thm-maxlcolsofACs}). 
(Nevertheless, Theorem~\ref{Thm-maxlcolsofACs} shows that \emph{most} Type~B cones only have a 1-dimensional intersection with $\mathcal{S}_0$.) 
The medium blue dot that is centered horizontally corresponds to the Type~2 curve of slope $-1$. 
It does not have an arc joining it to any green dot because there is no Type~1 curve of slope $-1$ whose shear coordinates sum to~$0$; however, it does have a gray arc joining it to the gray dot. 
This gray arc represents the 2-dimensional intersection of $\mathcal{S}_0$ with the Type~B cone whose Type~1 rays are given by $(-1, 1, -)$ and $(-1, 0, +)$. 
Finally, certain pairs of green dots are joined by green arcs. 
Theorem~\ref{Thm-compatible-arcs-distinct}, applied via the~$\kappa$ map to Type~1 curves given by either $(m, 1, -)$ with $m > -1$ or $(m, 0, +)$ with $m < -1$, implies that two green dots are joined by an arc if and only if the slopes $b/a$ and $d/c$ of their Type~1 curves satisfy $ad-bc = \pm 1$. 
Thus, three green dots form the vertices of a triangle if and only if the slopes $b/a$, $d/c$, and $f/e$ of their Type~1 curves satisfy $ad-bc = af-be = cf-de = 1$ (assuming $b/a < d/c < f/e$). 
Each of these triangles represents a 3-dimensional cone that is a face of exactly two Type~A cones (cf. Theorems~\ref{Thm-maxlcolsofACs} and ~\ref{thm-triangulations}). 
(However, Theorem~\ref{Thm-maxlcolsofACs} shows that \emph{most} Type~A cones do \emph{not} intersect $\mathcal{S}_0$.) 
Some green dots also have gray arcs joining them to the gray dot. 
Each of these gray arcs represents the 2-dimensional intersection of $\mathcal{S}_0$ and a 3-dimensional face shared by two Type~A cones that themselves intersect $\mathcal{S}_0$ in 3-dimensional cross-sections (each manifested as a triangle with two gray sides and one green side in both figures) on either side of the 2-dimensional intersection. 

Unfortunately, since $\mathcal{S}_1$ and $\mathcal{S}_{-1}$ are not linear, their intersections with the cones of $\mathcal{F}_\mathbb{Q}(\Delta_0)$ cannot be visualized on (2-dimensional) paper in a way that is both geometrically accurate and combinatorially meaningful. 
However, if we dispense with geometric accuracy, the entirety of $\mathcal{F}_\mathbb{Q}(\Delta_0)$ can be illustrated in a way that is combinatorially meaningful. 
For such a figure, see \cite{mydissertation}. 

Now that we have finished illustrating the geometry and combinatorics of $\mathcal{F}_{\mathbb{Q}}(\Delta_0)$, we make conjectures regarding (a) the nature of the non-rational cones of $\mathcal{F}_{B(\Delta_0)}$ and (b) which, if any, of the Type~C cones and/or Type~3 ray are in $\mathcal{F}_{B(\Delta_0)}$. 

\begin{conj} \label{conj-rational-fan-fills-NOTSUM=0}
Every point of $\mathbb{R}^4 \setminus \mathcal{S}_0$ is contained in a cone of $\mathcal{F}_{\mathbb{Q}}(\Delta_0)$.
\end{conj}
 
A possible approach for a proof of Conjecture~\ref{conj-rational-fan-fills-NOTSUM=0}  is sketched in \cite{mydissertation}. 
Before stating the next conjecture, we make a few observations. 
Let $C_{m_1, m_2}$ denote the nonnegative span of the shear coordinate vectors of the Type~2 curves of slopes $m_1$ and $m_2$. 
Let $C_\triangledown = C_{\infty,0} \cup C_{0,-1} \cup C_{-1,\infty}$. 
By Theorem~\ref{Thm-shear-coords-of-ACs}, the Type~2 rays lie in $C_\triangledown$ (it is the outline of the ``shield'' shape in Figure~\ref{SUM=0}) and comprise all the rational rays in $C_\triangledown$. 
Moreover, since the rational numbers are dense in the real numbers, Theorem~\ref{Thm-shear-coords-of-ACs} implies that the Type~2 rays are dense in $C_\triangledown$. 


\begin{conj} \label{conj-only-non-rational-cones-of-mut-fan}
The non-rational cones of $\mathcal{F}_{B(\Delta_0)}$ consist of:
\begin{enumerate}
\item \label{non-rational-rays} Each non-rational ray in $C_\triangledown$. 
\item \label{non-rational-2Dcones} For each non-rational ray in $C_\triangledown$, the nonnegative span of that ray and $[0, 0, 1, -1]$. 
\end{enumerate}
\end{conj}
 
Conjectures~\ref{conj-rational-fan-fills-NOTSUM=0} and~\ref{conj-only-non-rational-cones-of-mut-fan} are logically equivalent. 
This equivalence is proven in~\cite{mydissertation}. 

\begin{conj} \label{conj-all-in-mut}
All cones in $\mathcal{F}_{\mathbb{Q}(\Delta_0)}$ are in $\mathcal{F}_{B(\Delta_0)}$. 
\end{conj}

This last conjecture would follow as a corollary to the previous two (equivalent) conjectures. 
This implication is proven in \cite{mydissertation}. 
%

\end{document}